\documentstyle[12pt,amssymb,amscd]{amsart}  
\input amssym.def
\input amssym.tex

\newtheorem{thm}{Theorem}[subsection]

\newtheorem{cor}[thm]{Corollary}

\newtheorem{lemma}[thm]{Lemma}

\newtheorem{definition}[thm]{Definition}
\newtheorem{proposition}[thm]{Proposition}

\theoremstyle{remark}
\newtheorem{remark}[thm]{Remark}

\theoremstyle{definition}

\numberwithin{equation}{section}

\newcommand{\nc}{\nabla _{{\operatorname {can}}}}
\newcommand{\C}{{\Bbb C}}

\newcommand{\cO}{{\cal O}}

\newcommand{\cL}{{\cal L}}
\newcommand{\cM}{{\cal M}}

\newcommand{\cA}{{\cal A}}

\newcommand{\ea}{{\operatorname{exp}}\,{\operatorname{ad}}}
\newcommand{\tn}{{\widetilde{\nabla}}}

\newcommand{\cC}{{\cal C}}

\newcommand{\Tot}{{\operatorname{Tot}}}

\newcommand{\liminv}{{\operatorname{lim}}\,{\operatorname{inv}}}
\newcommand{\limdir}{{\operatorname{lim}}\,{\operatorname{dir}}}

\newcommand{\A}{{\cal{A}}}

\newcommand{\congr}{=}
\newcommand{\odsb}{\Omega_{\operatorname{DRS}}^{\bullet}}
\newcommand{\ods}{\Omega_{\operatorname{DRS}}}

\newcommand{\isomoto}{\overset{\sim}{\to}}

\newcommand{\g}{{\frak{g}}}

\newcommand{\mtwa}{{\operatorname{Matr}}_{\operatorname{tw}}(\cA)}

\newcommand{\no}{\nabla_0}

\date{July 13, 2005}

\begin{document}

{\centerline{\large{\bf DEFORMATION QUANTIZATION OF GERBES}}}
\smallskip
{\centerline{Paul Bressler  }}
\smallskip
{\centerline{Department of Mathematics, University of Arizona }}
\smallskip
{\centerline{ Alexander Gorokhovsky}}
\smallskip
{\centerline{Department of Mathematics, University of Colorado  }}
\smallskip
{\centerline{Ryszard Nest }}
\smallskip
{\centerline{ Mathematical Institute, University of Copenhagen}}
\smallskip
{\centerline{Boris Tsygan  }}
\smallskip
{\centerline{Department of Mathematics, Northwestern University}}
\bigskip
{\bf \centerline{Abstract}}

This is the first in a series of articles devoted to deformation quantization of gerbes. We introduce basic definitions, interpret deformations of a given stack as Maurer-Cartan elements of a differential graded Lie algebra (DGLA), and classify deformations of a given gerbe in terms of Maurer-Cartan elements of the DGLA of Hochschild cochains twisted by the cohomology class of the gerbe. We also classify all deformations of a given gerbe on a symplectic manifold, as well as provide a deformation-theoretic interpretation of the first Rozansky-Witten class.

\section{Introduction}  \label{s:introduction}

The notion of deformation quantization, as well as the term, was first introduced in \cite{BFFLS}. Both became standard since then. A deformation quantization of a manifold $M$ is a multiplication law on the ring of functions on $M$ which depends on a formal parameter $\hbar$. This multiplication law is supposed to satisfiy certain properties, in particular its value at $\hbar=0$ must be equal to the usual multiplication. A deformation quantization defines a Poisson structure on $M$; therefore it is natural to talk about deformation quantization of Poisson manifolds. In the case when $M$ is a symplectic manifold, deformation quantizations of $C^{\infty}(M)$ were classified up to isomorphism in \cite{DWL},  \cite{Fe}, \cite{D}. In the case of a complex manifold $M$ with a holomorphic symplectic form, deformation quantizations of the sheaf of algebras $\cO _M$ are rather difficult to study. They were classified, under additional cohomological assumptions, in \cite{NT} (Theorem \ref{thm:classification of sheaves} of the present paper; cf. also \cite{BK} for the algebraic case). If one moves away from symplectic to general Poisson manifolds, the problem becomes much more complicated. All deformation quantizations of $\cO _M$ were classified by Kontsevich in \cite{K1}. For the algebraic case, cf. \cite{Y}.

In this paper we start a program of studying deformation quantization of stacks and gerbes. Stacks are a natural generalization of sheaves of algebras. They appear in geometry, microlocal analysis and mathematical physics, cf. \cite{Gi}, \cite{Br}, \cite{DP}, \cite{Ka}, \cite{MMS}, \cite{MR1}, \cite{MR2}, \cite{PS}, and other works. 

The main results of this paper are as follows.

1) We prove that deformations of every stack (in the generality adopted by us here) are classified by Maurer-Cartan elements of a differential graded Lie algebra, or DGLA (Theorems \ref{thm:deformations of stacks via dgla}, \ref{thm:deformations of stacks via dgla II}).
This generalizes the results of Gerstenhaber \cite{Ge} for associative algebras and of Hinich \cite{H} for sheaves of associative algebras. 

2) We show that the DGLA controlling deformations of a gerbe on a manifold is equivalent to the Hochschild cochain complex of this manifold, twisted by the cohomology class of the gerbe (Theorem \ref{thm:classification of deformations of a gerbe}).

3) We classify deformation quantizations of all gerbes on a symplectic manifold (Theorems \ref{thm:classification of deformations of the trivial gerbe, symplectic case} and \ref{thm:symplectic classification}). This generalizes the classification results for deformation quantizations of $C^{\infty}$ symplectic manifolds \cite{DWL}, \cite{D}, \cite{Fe}, \cite{Fe1}.

4) We show that the first Rozansky-Witten class of a holomorphic symplectic manifold is an obstruction for a canonical stack deformation quantization to be a sheaf of algebras (Theorem \ref{thm:RW}).

We start by defining stacks, gerbes and their deformations in the generality suited for our purposes (section \ref{s:stacks and cocycles}). We then recall (in subsections \ref{ss:definitions of deformations}, \ref{ss:Descent data for Deligne two-groupoids} ) the language of differential graded Lie algebras (DGLAs) in deformation theory, along the lines of \cite{GM}, \cite{Ge}, \cite{S}, \cite{SS}, \cite{Dr}, \cite{HS}. Then we pass to a generality that suits us better, namely to the case of cosimplicial DGLAs (subsection \ref{ss:cosimplicial DGLAs}). We define descent data for the Deligne two-groupoid (cf. \cite{G}, \cite{G1} and references thereof) of a cosimplicial DGLA and prove that the set of isomorphism classes of such data does not change if one passes to a quasi-isomorphic cosimplicial DGLA (Proposition \ref{prop:quis of cosimplicial dglas}). Next, we recall the construction of totalization of a cosimplicial DGLA (subsection \ref{ss:Totalization of cosimplicial DGLAs}). We prove that isomorphism classes of descent data of a cosimplicial DGLA are in one-to-one correspondence with isomorphism classes of Maurer-Cartan elements of its totalization.

After that, given a gerbe on a Poisson manifold, we define its deformation quantization. We first classify deformations of the trivial gerbe, i.e. deformations of the structure sheaf as a stack, on a symplectic manifold $M$, $C^{\infty}$ or complex (Theorem \ref{thm:classification of deformations of the trivial gerbe, symplectic case}; this result is very close to the main theorem of \cite{P}). More precisely, we first reduce the classification problem to classifying certain {\em $Q$-algebras}, using the term of A.~Schwarz (or {\em curved DGAs}, as they are called in \cite{Bl}). (Similar objects were studied in several contexts, in particular in \cite{C}). The link between these objects and gerbes was rather well understood for some time; for example, it is through such objects that gerbes appear in \cite{Kapu}). We also give a new proof of the classification theorem for deformations of the sheaf of algebras of functions (Theorem \ref{thm:classification of sheaves}). Then we show how the first Rozansky-Witten class \cite{RW}, \cite{Kap}, \cite{K2}) can be interpreted as an obstruction for a certain canonical deformation of the trivial gerbe to be a sheaf, not just a stack. This canonical stack is very closely related to stacks of microdifferential operators defined in \cite{Ka} and \cite{PS}. 

Next, we show how to interpret deformations of any gerbe in the language of DGLAs (Theorems \ref{thm:classification of deformations of a gerbe} and \ref{thm:classification of deformations of a holomorphic gerbe}). The proof is based on a DGLA interpretation of the deformation theory of any stack (within our generality); such an interpretation is provided by Theorem \ref{thm:deformations of stacks via dgla}. We show there that deformations of a stack are classified by the DGLA of De Rham-Sullivan forms with coefficients in {\em local Hochschild cochains of the twisted matrix algebra} associated to this stack. 

Note that De Rham-Sullivan forms were used in \cite{Y} to classify deformation quantizations of algebraic varieties. 

(The DGLA above is actually a DGLA of Hochschild cochains of a special kind of an associative DGA; the cyclic homology of this DGA is the natural recipient of the Chern character of a twisted module over a stack. We will study this in the sequel). 

Afterwards we prove a classification theorem for deformation quantizations of any gerbe on a symplectic manifold (Theorems \ref{thm:symplectic classification} and \ref{thm:symplectic holomorphic classification}). This can be viewed as an adaptation of Fedosov's methods \cite{Fe}, \cite{Fe1} to the case of gerbes. Note that some ideas about deformation quantization of gerbes appeared already in Fedosov's work; cf. also \cite{K}, as well as \cite{Ka} and \cite{PS}.

This paper was motivated by the index theory, in particular by index theorems for Fourier integral operators or by index theorems such as in \cite{MMS}. Among the applications other than index theory, we would like to mention dualities between gerbes and noncommutative spaces, as in \cite{Kapu}, \cite{Bl}, \cite{BBP}, \cite{MR1}, \cite{MR2}. The deformation-theoretical role of the first Rozansky-Witten class is also quite intriguing and worthy of further study.

The research of A.~G. and B.~T. was partially supported by NSF grants.
\section{Stacks and cocycles} \label{s:stacks and cocycles}
\subsection{}
Let $M$ be a smooth manifold ($C^{\infty}$ or complex). In this paper, by a stack on $M$ we will mean the following data:

1) an open cover $M=\cup U_i$;

2) a sheaf of rings $\A _i$ on every $U_i$;

3) an isomorphism of sheaves of rings $G_{ij}: \A _j|(U_i \cap U_j) \isomoto \A _i |(U_i \cap U_j)$ for every $i,\,j$;

4) an invertible element $c_{ijk} \in \A _i (U_i \cap U_j \cap U_k)$ for every $i,\,j,\,k$ satisfying 
\begin{equation} \label {eq:2-cocycle 1}
G_{ij}G_{jk}={\operatorname{ Ad}}(c_{ijk})G_{ik}
\end{equation}
such that, for every $i,\,j,\,k,\,l$,
\begin{equation} \label {eq:2-cocycle 2}
c_{ijk}c_{ikl}=G_{ij}(c_{jkl})c_{ijl}
\end{equation}

If two such data $(U'_i,\; \A' _i, \; G'_{ij},\; c'_{ijk})$ and $(U''_i,\; \A'' _i, \; G''_{ij},\; c''_{ijk})$ are given on $M$, an isomorphism between them is an open cover $M=\cup U_i$ refining both $\{U'_i\}$ and $\{U''_i\}$ together with isomorphisms $H_i: \A' _i \isomoto \A''_i$ on $U_i$ and invertible elements $b_{ij}$ of $\A' _i (U_i \cap U_j)$ such that
\begin{equation} \label{eq:equivalence of stacks 1}
G''_{ij}=H_i {\operatorname{Ad}}(b_{ij})G'_{ij}H_j ^{-1}
\end{equation}
and
 \begin{equation} \label{eq:equivalence of stacks 2}
H_i^{-1}(c''_{ijk})=b_{ij}G'_{ij}(b_{jk})c'_{ijk}b_{ik}^{-1}
\end{equation}

{\em A gerbe} is a stack for which $\A _i = {\cal O}_{U_i}$ and $G_{ij}={\operatorname {id}}$. In this case $c_{ijk}$ form a two-cocycle in $Z^2(M,{\cal O}_M^*)$.
\subsection{Categorical interpretation}\label{ss:categorical} Here we remind the well-known categorical interpretation of the notions introduced above. Though not used in the rest of the paper, this interpretation provides a very strong motivation for what follows.

A stack defined as above gives rise to the following categorical data:

1) A sheaf of categories $\cC _i$ on $U_i$ for every $i$;

2) an invertible functor $G_{ij}: \cC _j |(U_i \cap U_j) \isomoto \cC _i |(U_i \cap U_j)$ for every $i,\,j$;

3) an invertible natural transformation 
$$c_{ijk}: G_{ij}G_{jk}|(U_i \cap U_j \cap U_k) \isomoto G_{ik}|(U_i \cap U_j \cap U_k)$$
such that, for any $i,\,j,\,k,\,l$, the two natural transformations from 
\newline $G_{ij}G_{jk}G_{kl}$ to $G_{il}$ that one can obtain from the $c_{ijk}$'s are the same on $U_i \cap U_j \cap U_k \cap U_l$.

If two such categorical data $(U'_i,\; \cC' _i, \; G'_{ij},\; c'_{ijk})$ and $(U''_i,\; \cC'' _i, \; G''_{ij},\; c''_{ijk})$ are given on $M$, an isomorphism between them is an open cover $M=\cup U_i$ refining both $\{U'_i\}$ and $\{U''_i\}$, together with invertible functors $H_i: \cC' _i \isomoto \cC''_i$ on $U_i$ and invertible natural transformations $b_{ij}:H_i G'_{ij}|(U_i \cap U_j) \isomoto G''_{ij}H_j|(U_i \cap U_j)$ such that, on any $U_i \cap U_j \cap U_k$, the two natural transformations $H_i G'_{ij}G'_{jk} \isomoto G''_{ij}G''_{jk}H_k$ that can be obtained using $H_i$'s, $ b_{ij}$'s, and $c_{ijk}$'s are the same. More precisely:
 \begin{equation} \label{eq:equivalence of categorical data for stacks 1}
((c''_{ijk})^{-1}H_k)(b_{ik})(H_i c'_{ijk})=(G''_{ij}b_{jk})(b_{ij} G'_{jk})
\end{equation}

The above categorical data are defined from $(\A _i, G_{ij}, c_{ijk})$ as follows:

1) $\cC _i$ is the sheaf of categories of $\A _i$-modules;

2) given an $\A _i$-module ${\cal M}$, the $\A _j$-module $G_{ij}({\cal M})$ is the sheaf ${\cal M}$ on which $a \in \A _i$ acts via $G_{ij}^{-1}(a)$;

3) the natural transformation $c_{ijk}$ between $G_{ij}G_{jk}({\cal M})$ and $G_{jk}({\cal M})$ is given by multiplication by $G_{ik}^{-1}(c_{ijk}^{-1})$.

From the categorical data defined above, one defines a sheaf of categories on $M$ as follows. For an open $V$ in $M$, an object of $\cC (V)$ is a collection of objects $X_i$ of $\cC _i (U_i \cap V)$, together with isomorphisms $g_{ij}: G_{ij}(X_j) \isomoto X_i$ on every $U_i \cap U_j \cap V$, such that
$$g_{ij}G_{ij}(g_{jk})=g_{ik}c_{ijk}$$ on every $U_i \cap U_j \cap U_k \cap V$. A morphism between objects $(X'_i, g'_{ij})$ and $(X''_i, g''_{ij})$ is a collection of morphisms $f_i: X'_i \to X''_i$ (defined for some common refinement of the covers), such that $f_i g'_{ij}=g''_{ij}G_{ij}(f_j)$.

\begin{remark}\label{rmk:terminology of stacks} What we call stacks are what is referred to in \cite{DP} as descent data for a special kind of stacks of twisted modules (cf. Remark 1.9 in \cite{DP}). Both gerbes and their deformations are stacks of this special kind. We hope that our terminology, which blurs the distinction between stacks and their descent data, will not cause any confusion.
\end{remark}
\subsection{Deformations of stacks}
\begin{definition} \label{dfn:deformation1}
Let $k$ be a field of characteristic zero. Let ${\mathfrak{a}}$ be a local Artinian k-algebra with the maximal ideal ${\mathfrak{m}}$. A deformation of a stack $\cA ^{(0)}$ over ${\mathfrak{a}}$ is a stack $\cA$ where all $\cA _i$ are sheaves of ${\mathfrak{a}}$-algebras, free as ${\mathfrak{a}}$-modules, $G_{ij}$ are isomorphisms of algebras over ${\mathfrak{a}}$, and the induced stack $\cA/{\mathfrak{m}}\cA$ is equal to $\cA ^{(0)}$. An isomorphism of two deformations is an isomorphism of stacks which is identity modulo ${\mathfrak{m}}$ and such that $H_i$ are isomorphisms of algebras over ${\mathfrak{a}}.$
\end{definition}
Consider the filtration of ${\mathfrak{a}}$ by powers of ${\mathfrak{m}}.$ Choose a splitting of the filtered $k$-vector space
$${\mathfrak{a}}=\oplus_{m=0}^N {\mathfrak{m}}_m$$
where ${\mathfrak{m}}_m={\mathfrak{m}}^m/{\mathfrak{m}}^{m+1}.$ 

Given a deformation, we can identify $\cA_i=\cA^{(0)}_i\otimes {\mathfrak{a}};$ the multiplication on $\cA_i$ is determined by
$$f*_i g = fg + \sum_{m=1}^N  P_i^{(m)}(f,g)$$
with $P_i^{(m)}: {\cA^{(0)}_i}^{\otimes 2}\to \cA^{(0)}_i\otimes {\mathfrak{m}}_m.$ Similarly, $G_{ij}$ is determined by 
$$G_{ij}(f) =  f + \sum_{m=1}^{N}  T_{ij}^{(m)}(f)$$
with $T_{ij}^{(m)}:\cA^{(0)}\to \cA^{(0)}_i\otimes {\mathfrak{m}}_m,$ and 
$$c_{ijk} = \sum _{m=0}^{N} c^{(m)}_{ijk}$$
with $c^{(m)}_{ijk}\in {\mathfrak{m}}_m.$
For an isomorphism of two stacks, $H_i$ is determined by 
$$H_{i}(f) =  f + \sum_{m=1}^{N}  R_{i}^{(m)}(f)$$
with $H_{i}^{(m)}:\cA^{(0)}\to \cA^{(0)}_i\otimes {\mathfrak{m}}_m;$
$$b_{ij} = \sum _{m=0}^{N} b^{(m)}_{ij}$$
with $b^{(m)}_{ij}\in {\mathfrak{m}}_m.$
\begin{definition} \label{dfn:deformation}
Consider a gerbe $\cA^{(0)}$ given by a two-cocycle $c^{(0)}_{ijk}$. A deformation of $\cA^{(0)}$ is by definition its deformation as a stack, such that $P_i^{(m)}(f,g)$ are (holomorphic) bidifferential expressions and $T_{ij}^{(m)}$are (holomorphic) differential expressions.

An isomorphism between two deformations is an isomorphism $(H_i, b_{ij})$
where $R_i^{(m)}$ are (holomorphic) differential expressions.
\end{definition}

\section{Differential graded Lie algebras and deformations} \label{s:Differential graded Lie algebras and deformations}
\subsection{} \label{ss:definitions of deformations}
Here we give some definitions that lie at the foundation of the deformation theory program along the lines of \cite{Ge}, \cite{GM}, \cite{S}, \cite{SS}, \cite{Dr}, \cite{HS}, as well as of the notions such as Deligne two-groupoid (cf. \cite{G}, \cite{G1} and references thereof). Let
$$\cL = \bigoplus_{m\geq -1} \cL ^m$$
be a differential graded Lie algebra (DGLA). Let ${\mathfrak{a}}$ be a local Artinian k-algebra with the maximal ideal ${\mathfrak{m}}$. We call {\em a Maurer-Cartan element}  an element ${\lambda}$ of $\cL ^1 \otimes {\mathfrak{m}}$ satisfying 
\begin{equation} \label{eq:MC}
d\lambda + \frac{1}{2}[\lambda,\lambda]=0
\end{equation}
{\em A gauge equivalence} between two Maurer-Cartan elements $\lambda $ and $\mu$ is an element $G={\operatorname {exp}}\,X$ where $X\in \cL ^0 \otimes {\mathfrak{m}}$ such that
\begin{equation} \label{eq:MC equivalence}
d+\mu={\operatorname {exp}}\,{\operatorname {ad}}X \,(d+\lambda)
\end{equation}
The latter equality takes place in the cross product of the one-dimensional graded Lie algebra $kd$ concentrated in dimension one and $\cL ^0 \otimes {\mathfrak{m}}$.
Given two gauge transformations $G={\operatorname{exp}}\,X,\,H={\operatorname{exp}}\, Y$ between $\lambda$ and $\mu$, {\em a two-morphism} from $H$ to $G$ is an element $c= {\operatorname{exp}}\,t$ of $\cL ^{-1} \otimes{\mathfrak{m}} $ such that 
\begin{equation} \label{eq:MC 2-morphism}
{\operatorname {exp}}(X)={\operatorname {exp}}(dt+[\mu,t]) {\operatorname {exp}}Y
\end{equation}
in the unipotent group ${\operatorname {exp}} (\cL ^0\otimes {\mathfrak{m}})$. The composition of gauge transformations $G$ and $H$ is the product $GH$ in the unipotent group ${\operatorname {exp}} (\cL ^0)\otimes {\mathfrak{m}}$. The composition of two-morphisms $c_1$ and $c_2$ is the product $c_1c_2$ in the prounipotent group ${\operatorname {exp}}(\cL ^{-1}\otimes {\mathfrak{m}})$. Here $\cL ^{-1}\otimes {\mathfrak{m}}$ is viewed as a Lie algebra with the bracket 
\begin{equation}\label{eq:mu-bracket}
[a,\,b]_{\mu}=[a,\,\delta b+[\mu, \,b]]
\end{equation}
We denote the above pronilpotent Lie algebra by $(\cL ^{-1}\otimes {\mathfrak{m}})_{\mu}.$
 The above definitions, together with the composition, provide the definition of {\em the Deligne two-groupoid} of $\cL\otimes {\mathfrak{m}}$ (cf. \cite{G1})..

\begin{remark}\label{rmk:ezra}
Recently Getzler gave a definition of a Deligne $n$-groupoid of a DGLA concentrated in degrees above $-n$, cf. \cite{G}.
\end{remark}

\subsection{Descent data for Deligne two-groupoids} \label{ss:Descent data for Deligne two-groupoids} Let $\cL$ be a sheaf of DGLAs on $M$. A descent datum of the Deligne two-groupoid of $\cL\otimes {\mathfrak{m}}$ are the following: 

1) A Maurer-Cartan element $\lambda _i \in \cL^{1} \otimes {\mathfrak{m}}$ on $U_i$ for every $i$;

2) a gauge transformation $G_{ij}: \lambda_j |(U_i \cap U_j) \isomoto \lambda_i |(U_i \cap U_j)$ for every $i,\,j$;

3) a two-morphism
$$c_{ijk}: G_{ij}G_{jk}|(U_i \cap U_j \cap U_k) \isomoto G_{ik}|(U_i \cap U_j \cap U_k)$$
such that, for any $i,\,j,\,k,\,l$, the two two-morphisms from $G_{ij}G_{jk}G_{kl}$ to $G_{il}$ that one can obtain from the $c_{ijk}$'s are the same on $U_i \cap U_j \cap U_k \cap U_l$.

If two such data $(U'_i,\; \lambda' _i, \; G'_{ij},\; c'_{ijk})$ and $(U''_i,\; \lambda'' _i, \; G''_{ij},\; c''_{ijk})$ are given on $M$, an isomorphism between them is an open cover $M=\cup U_i$ refining both $\{U'_i\}$ and $\{U''_i\}$, together with gauge transformations $H_i: \lambda' _i \isomoto \lambda''_i$ on $U_i$ and two-morphisms $b_{ij}:H_i G'_{ij}|(U_i \cap U_j) \isomoto G''_{ij}H_j|(U_i \cap U_j)$ such that, on any $U_i \cap U_j \cap U_k$, the two two-morphisms $H_i G'_{ij}G'_{jk} \isomoto G''_{ij}G''_{jk}H_k$ that can be obtained using $H_i$'s, $ b_{ij}$'s, and $c_{ijk}$'s are the same.

Finally, given two isomorphisms $(H'_i, b'_{ij})$ and $(H''_i, b''_{ij})$ between the two data $(U_i,\; \lambda' _i, \; G'_{ij},\; c'_{ijk})$ and $(U_i,\; \lambda'' _i, \; G''_{ij},\; c''_{ijk})$, define a {\em two-isomorphism} between them to be a collection of two-morphisms $a_i: H'_i \to H''_i$ such that 
$$b''_{ij}\circ (a_i \circ G'_{ij})=(G''_{ij}\circ a_i)\circ b'_{ij}$$
as two-morphisms from $H'_i\circ G'_{ij}\to G''_{ij}\circ H''_{j}$.
\subsection{Cosimplicial DGLAs and descent data} \label{ss:cosimplicial DGLAs}
The notion of a descent datum above, as well as an analogous notion for simplicial sheaves of DGLAs that we use below, is a partial case of a more general situation that we are about to discuss. Recall that {\em a cosimplicial object} of a category $\cC$ is a functor $X: \Delta \to \cC$ where $\Delta$ is the category whose objects are sets $[n]=\{0,\,\ldots , \, n\}$ with the standard linear ordering ($n\geq 0$), and morphisms are nondecreasing maps. We denote $X([n])$ by $X^n.$ For $0\leq i\leq n,$ let $d_i:[n]\to [n+1]$ be the only injective map such that $i$ is not in the image, and $s_i: [n+1]\to [n]$ the only surjection for which every element of $[n-1]$ except $i$ has exactly one preimage. For a cosimplicial Abelian group $\cA$, one defines the standard differential
$$\partial=\sum_{i=0}^n(-1)^id_i:\cA^n\to \cA^{n+1}.$$

For a cosimplicial set $X,$ let $x\in X^k.$ Let $n\geq k$ and $0\leq i_0<\ldots <i_k\leq n.$ By $x_{i_0\ldots i_k}$ we denote the object of $X^n$ which is the image of $x$ under the map in $\Delta$ which embeds $[k]$ into $[n]$ as the subset $\{i_0,\,\ldots,\, i_k\}.$

Let $\cL$ be a cosimplicial DGLA. We will denote by $\cL ^{n,p}$ the component of degree $p$ of the DGLA $\cL ^n$, $n\geq 0.$ 

Let ${\mathfrak{a}}$ be a local Artinian algebra over $k$ with the maxiamal ideal ${\mathfrak{m}}$. Consider a cosimplicial DGLA $\cL$ such that $\cL ^{n,p}=0$ for $p<-1.$ {\em A descent datum for the Deligne two-groupoid of  $\cL\otimes {\mathfrak{m}}$} is the following:

1) A Maurer-Cartan element $\lambda  \in \cL^{0,1} \otimes {\mathfrak{m}};$

2) a gauge transformation $G: \lambda_1  \isomoto \lambda_0$ in ${\operatorname{exp}}(\cL^{1,0});$

3) a two-morphism
$$c: G_{01}G_{12} \isomoto G_{02}$$
in ${\operatorname{exp}}(\cL^{2,-1}_{\lambda_0})$
such that, for any $i,\,j,\,k,\,l$, the two two-morphisms from $G_{01}G_{12}G_{23}$ to $G_{03}$ that one can obtain from the $c_{ijk}$'s are the same.

An isomorphism between two data $(\lambda' , \; G',\; c')$ and $( \lambda'' , \; G'',\; c'')$ is a pair of a gauge transformation $H: \lambda'  \isomoto \lambda''$ and a two-morphism $b_{01}:H_0 G'_{01} \isomoto G''_{01}H_1$ such that the two two-morphisms $H_0 G'_{01}G'_{12} \isomoto G''_{01}G''_{12}H_2$ that can be obtained using $H_i$'s, $ b_{ij}$'s, and $c$ are the same.

For two isomorphisms $(H', b')$ and $(H'', b'')$ between the two data $ \lambda' , \; G',\; c')$ and $(U,\; \lambda'' , \; G'',\; c'')$, define a {\em two-isomorphism} between them to be a collection of two-morphisms $a: H' \to H''$ such that 
$$b''_{01}\circ (a_0 \circ G'_{01})=(G''_{01}\circ a_0)\circ b'_{01}$$
as two-morphisms from $H'_0\circ G'_{01}\to G''_{01}\circ H''_{1}$.

\begin{proposition} \label{prop:quis of cosimplicial dglas}
a). A morphism $f:\cL _1 \to \cL _2$ of cosimplicial DGLAs induces a map from the set of isomorphism classes of descent data of the Deligne two-groupoid of $\cL _1\otimes {\mathfrak{m}}$ to the set of isomorphism classes of descent data of the Deligne two-groupoid of $\cL _2\otimes {\mathfrak{m}}$. 

b). Assume that $f$ induces a quasi-isomorphism of total complexes of the double complexes $\cL _1^{n,p} \to \cL _2^{n,p} $. Then the map defined in a) is a bijection.

c). Under the assumptions of b), let $\cA$ be a descent datum of the Deligne two-groupoid of $\cL _1\otimes {\mathfrak{m}}$, and let $f(\cA)$ be its image under the map from a). The morphism $f$ induces a bijection
$$\frac{{\operatorname{Iso}}(\cA, \cA')}{2-{\operatorname{Iso}}}\isomoto \frac{{\operatorname{Iso}}(f(\cA), f(\cA'))}{2-{\operatorname{Iso}}}.$$

d). For two isomorphisms $\phi, \psi:\cA\to \cA'$, denote their images under the above bijection by $f(\phi), f(\psi).$ Then $f$ induces a bijection
$$2-{\operatorname{Iso}}(\phi, \psi)\isomoto 2-{\operatorname{Iso}}(f(\phi), f(\psi))$$
\end{proposition}
In other words, $f$ induces {\em an equivalence of two-groupoids of descent data}, compare to \cite{G}, \cite{G1}.

{\bf Proof.} What follows is essentially a standard deformation theoretical proof. We start by establishing a rigorous expression of the following intuitive statement. First, a descent datum $(\lambda, G, c)$ is a non-Abelian version of a two-cocycle of the double complex $\cL^{\bullet, \bullet}, \partial+d;$ second, if one takes an arbitrary datum $(\lambda, G, c)$ and measures its deviation from being a descent datum, the result will be a non-Abelian version of a three-cocycle. This is, in essence, what enables us to study deformations of descenta data by homological methods.

\subsubsection{}\label{sss:a}Consider a triple $(\lambda, G, t)$ with $\lambda \in \cL^{0,1}\otimes {\mathfrak m},$ $G \in {\operatorname{exp}}(\cL^{1,0}\otimes {\mathfrak m}),$ and $t \in \cL^{2,-1}\otimes {\mathfrak m}.$ Define the operation $a\cdot _{\lambda} b$ on $\cL^{-1}\otimes {\mathfrak m}$ to be the Campbell-Dynkin-Hausdorff series corresponding to the bracket $[a,\,b]_{\lambda}=[a,\,db+[\lambda,b]].$ If $\lambda$ is a Maurer-Cartan element, this is a group multiplication. If not, one can still define the operation which is no longer associative; zero is the neutral element, and every element is invertible. Denote the set $\cL^{-1}\otimes {\mathfrak m}$ with the operation $a\cdot _{\lambda} b$ by ${\operatorname{exp}}((\cL^{-1}\otimes {\mathfrak m})_{\lambda})$. For $t \in \cL^{-1}\otimes {\mathfrak m},$ we will denote by ${\operatorname{exp}}(t)$ the element $t$ viewed as an element of ${\operatorname{exp}}((\cL^{-1}\otimes {\mathfrak m})_{\lambda})$.

{\bf A notation convention.} For $G\in {\operatorname{exp}}(\cL^{1,0}\otimes {\mathfrak m})$ and $X\in \cL\otimes {\mathfrak m}),$ we will denote ${\operatorname{Ad}}_G(X)$ simply by $G(X).$ For $\lambda \in \cL^{0,1}\otimes {\mathfrak m},$, $G(d+\lambda)$ will stand for $d+\lambda'$ where $\lambda'$ is the image of $\lambda$ under the gauge transformation by $\lambda.$

Given $(\lambda, G, t)$ as above; let $c={\operatorname{exp}}(t)$ and $\gamma={\operatorname{exp}}(dt+[\lambda_0,t])$ in ${\operatorname{exp}}(\cL^{1,0}\otimes {\mathfrak m}).$ Define

\begin{equation}\label{eq:R,Z,T,Phi 1}
R=d\lambda+\frac{1}{2}[\lambda,\lambda]\in \cL^{0,2}\otimes {\mathfrak m};
\end{equation}
\begin{equation}\label{eq:R,Z,T,Phi 2}
Z=G(d+\lambda_1)-(d+\lambda_0) \in \cL^{1,1}\otimes {\mathfrak m};
\end{equation}
\begin{equation}\label{eq:R,Z,T,Phi 3}
G_{02}=T\gamma G_{01} G_{12}\in {\operatorname{exp}}(\cL^{2,0}\otimes {\mathfrak m})
\end{equation}
(this is a definition of $T$);
\begin{equation}\label{eq:R,Z,T,Phi 4}
\Phi=((G_{01}(c_{123})^{-1}c_{013}^{-1})c_{023})c_{012}\in {\operatorname{exp}}((\cL^{3,-1}\otimes {\mathfrak m})_{\lambda_0})
\end{equation}
(the order of parentheses is in fact irrelevant for our purposes).

Define ${\mathcal I}$ to be the cosimplicial ideal of $\cL\otimes{\mathfrak m}$ generated by $[R_{i}, \cL\otimes{\mathfrak m}]$, $[Z_{ij}, \cL\otimes{\mathfrak m}]$, $({\operatorname{Ad}}(T_{ijk})-{\operatorname{Id}}) (\cL\otimes{\mathfrak m})$. Note that the operation $a\cdot _{\lambda} b$ becomes a group law modulo ${\operatorname{exp}}({\mathcal{I}})$.
\begin{lemma}\label{lemma:nonabelian d2=0}
1) (The Bianchi identity): $dR+[\lambda, R]=0;$

2) (Gauge invariance of the curvature): 
$$R_0+dZ+[\lambda_0, Z]+\frac{1}{2}[Z,Z]-G(R_1)=0;$$

3) $T\gamma(d+\lambda_0)-(d+\lambda_0)+Z_{01}+G_{01}(Z_{12})-Z_{02}=0;$

4) $T_{013}(\gamma_{013}G_{01})(T_{123})\gamma_{013}G_{01}(\gamma_{123})=T_{023}\gamma_{023}(T_{012})\gamma_{023}\gamma_{012}$ 

modulo ${\operatorname{exp}}({\mathcal{I}})$;

5) (The pentagon equation):
$$G_{01}(\Phi _{1234}){{\operatorname{Ad}}_{G_{01}(c_{123})^{-1}}}(\Phi _{0134})\Phi _{0123}={{\operatorname{Ad}}_{G_{01}G_{12}(c_{234})^{-1}}}(\Phi _{0124})\Phi _{0234}$$
modulo ${\operatorname{exp}}({\mathcal{I}})$.
\end{lemma}
{\bf Proof.} The first equality is straightforward. The second follows from $(G(d+\lambda_1))^2=G((d+\lambda_1)^2).$ The third is obtained by applying both sides of \eqref{eq:R,Z,T,Phi 3} to $d+\lambda_2.$ The fourth can be seen by transforming $G_{01}(G_{12}G_{23})=(G_{01}G_{12})G_{23}$ in two different orders, using \eqref{eq:R,Z,T,Phi 3}. The fifth equation compares two two-morphisms from $((G_{01}G_{12})G_{23})G_{34}$ to $G_{01}(G_{12}(G_{23}G_{34}))$ corresponding to the two different routes along the perimeter of the Stasheff pentagon. (This is just a motivation for writing the formula which is then checked directly. We could not think of a reason for this formula to be true {\em a priori}).
\begin{cor}\label{cor:non-abelian d2=0, II} Let $(\lambda, \,G, t)$ be as in the beginning of \ref{sss:a}. Assume that they define a descent datum modulo ${\mathfrak{m}}^{n+1}.$ Then $(R^{(n+1)},$ $ Z^{(n+1)},$ $ T^{(n+1)},$ $ -\Phi^{(n+1)})$ is a $d+\partial$-cocycle of degree three.
\end{cor} 
\subsubsection{}\label{sss:aa} We need analogues of the above statements for isomorphisms and two-morphisms. Let $(\lambda, G, c)$ and $(\lambda', G', c')$ be two descent data. Consider a pair $(H, s)$ where $H\in {\operatorname{exp}}(\cL^{0,0}\otimes {\mathfrak{m}})$ and $s\in \cL^{1,-1}\otimes {\mathfrak{m}}.$ Put $b={\operatorname{exp}}(s)$ in ${\operatorname{exp}}((\cL^{1,-1}\otimes {\mathfrak{m}})_{\lambda'_0})$. Define also $\beta={\operatorname{exp}}((ds+[{\lambda'}_0,s])$ in ${\operatorname{exp}}((\cL^{1,0}\otimes {\mathfrak{m}}))$.

As above, we measure the deviation of the pair $(H,b)$ from being an isomorphism of descent data. Put 
\begin{equation}\label{eq:C,S,Psi 1}
C=H(d+\lambda)-(d+{\lambda'}) \in \cL^{0,1}\otimes {\mathfrak m};
\end{equation}
\begin{equation}\label{eq:C,S,Psi 2}
H_0G=S\beta G'H_1\in {\operatorname{exp}}(\cL^{1,0}\otimes {\mathfrak m})
\end{equation}
(this is a definition of $S$);
\begin{equation}\label{eq:C,S,Psi 3}
\Psi=b_{02}^{-1}{c'}_{012}{G'}_{01}(b_{12})b_{01}H_0(c_{012}) 
\end{equation}
in ${\operatorname{exp}}((\cL^{2,-1}\otimes {\mathfrak m})_{\lambda'_0}).$ The pair $(H,b)$ is an isomorphism between the two descent data if and only if $C=0, $ $S=1,$ $\Psi=1.$

Denote by $\mathcal J$ the cosimplicial ideal of $\cL\otimes{\mathfrak m}$ generated by $[C_i,\cL\otimes{\mathfrak m}]$ and $({\operatorname{Ad}}(S_{ij})-{\operatorname{Id}}) (\cL\otimes{\mathfrak m}).$
\begin{lemma}\label{lemma:nonabelian d2=0 I}
1) $dC+[\lambda', C]+\frac{1}{2}[C,C]=0;$

2) $S\beta (d+\lambda'_0)-(d+\lambda'_0)+C_0-G'_{01}(C_1)=0;$

3) $S_{01}\beta_{01}G'_{01}(S_{12}\beta_{12})=H_0(\gamma_{012})S_{02}\beta_{02}{\gamma'}_{012}^{-1}$
modulo ${\operatorname{exp}}({\mathcal{J}})$;

4) $\Psi_{023}{\operatorname{Ad}}_{H_0(c_{023})}(\Psi_{012})={\operatorname{Ad}}_{{b_{03}}^{-1}{c_{013}}G'_{01}(b_{13})}(G'_{01}(\Psi_{123}))\Psi_{012}$

modulo ${\operatorname{exp}}({\mathcal{J}})$.
\end{lemma}
{\bf Proof.} The first equality follows from $(H(d+\lambda))^2=0;$ the second from comparing the action of both sides of \eqref{eq:C,S,Psi 2} on $d+\lambda'_1;$ the third is obtained by comparing two different expressions for $H_0G_{01}G_{12}$ that can be obtained from \eqref{eq:C,S,Psi 2}. The fourth equality compares two different two-morphisms from $H_0G_{03}$ to itself. If one passes to two-morphisms from $H_0G_{01}G_{12}G_{03}$ to itself, it becomes the pentagon equation which compares two different routes from $((H_0G_{01})G_{12})G_{03}$ to $H_0(G_{01}(G_{12}G_{03}))$. One side of the pentagon, namely the edge between $H_0((G_{01}G_{12})G_{03})$ and $H_0(G_{01}(G_{12}G_{03}))$, degenerates into a point.
\begin{cor}\label{cor:non-abelian d2=0, III} Let $(H, s)$ be as in the beginning of \ref{sss:aa}. Assume that they define an isomorphism of descent data $(\lambda, G, c)$ and $(\lambda', G', c')$ modulo ${\mathfrak{m}}^{n+1}.$ Then $(C^{(n+1)},$ $ S^{(n+1)},$ $-\Psi^{(n+1)})$ is a $d+\partial$-cocycle of degree two.
\end{cor} 
\subsubsection{} \label{sss:aaa}Finally, we need an analogous statement for two-morphisms. Let $(H,b)$ and $({\widetilde{H}}, {\widetilde{b}})$  be isomorphisms between the descent data $(\lambda, G, c)$ and $(\lambda', G', c')$. Let $r\in \cL^{0,{-1}}\otimes {\mathfrak m}$ and $a={\operatorname{exp}}(r)$ in ${\operatorname{exp}}((\cL^{0,-1}\otimes{\mathfrak m})_{\lambda'}).$ Define $P$ and $\Omega$ by
\begin{equation}\label{eq:P,Omega 1}
{\widetilde{H}}=P\alpha H
\end{equation}
where $\alpha={\operatorname{exp}}((d+\lambda')a).$ Let ${\mathcal{K}}$ be the cosimplicial ideal generated by all $({\operatorname{Ad}}_{P_i}-{\operatorname{Id}})(\cL\otimes {\mathfrak m}).$
\begin{equation}\label{eq:P,Omega 2}
{\widetilde{b}}_{01}a_0=\Omega G'_{01}(a_1)b_{01}.
\end{equation}
$a: (H,b)\to ({\widetilde{H}}, {\widetilde{b}})$ is a two-morphism if and only if $P=1$ and $\Omega=1.$
\begin{lemma}\label{lemma:nonabelian d2=0 II}
1) $(P\alpha)(d+\lambda')=0;$

2) ${\operatorname{Ad}}_{\widetilde{\beta}}(P_0^{-1})G'(P_1)=({\widetilde{\beta}}\alpha_0)(G'(\alpha_1)\beta)^{-1}$
where $\beta={\operatorname{exp}}((d+\lambda')b)$ and ${\widetilde{\beta}}={\operatorname{exp}}((d+\lambda'){\widetilde{b}});$

3) ${\operatorname{Ad}}_{G'_{01}({\widetilde{b}}_{12})}(\Omega_{01})G'_{01}(\Omega_{12}){\operatorname{Ad}}_{c'_{012}}(\Omega_{02}^{-1})=1$
modulo ${\operatorname{exp}}({\mathcal{K}})$.
\end{lemma}
{\bf Proof.} The first equality follows from the fact that $H$ and ${\widetilde{H}}$ both preserve $d+\lambda'.$ The second is obtained by compairing the equalities $H_0G=\beta G'H_1,$ ${\widetilde{H}}_0G={\widetilde{\beta}}G'{\widetilde{H}}_1,$ and \eqref{eq:P,Omega 1}. The third equality is obtained by comparing two different expressions for $G'_{01}({\widetilde{b}}_{12}){\widetilde{b}}_{01}a_0$ using \eqref{eq:P,Omega 2}.
\begin{cor}\label{cor:non-abelian d2=0, IV} Let $r$ be as in the beginning of \ref{sss:aaa}. Assume that it defines a two-isomorphism $(H,b)\to ({\widetilde{H}}, {\widetilde{b}})$ modulo ${\mathfrak{m}}^{n+1}.$ Then $(P^{(n+1)},$ $-\Omega^{(n+1)})$ is a $d+\partial$-cocycle of degree one.
\end{cor} 
\subsubsection{End of the proof of Proposition \ref{prop:quis of cosimplicial dglas}} The statement a) is obvious. Let us prove the surjectivity of b). Let $(\mu,$ $G,$ $c$ be a descent datum for $\cL_2.$ Let $G={\operatorname{exp}}(y)$ in ${\operatorname{exp}}(\cL^{1,-0}_2\otimes{\mathfrak m})$ and $c={\operatorname{exp}}(t)$ in ${\operatorname{exp}}((\cL^{2,-1}_2\otimes{\mathfrak m})_\mu).$ We write
$$y=\sum y^{(k)};\;t=\sum t^{(k)}$$
{\em etc.}, where $y^{(k)},\,t^{(k)}\in \cL_2\otimes{\mathfrak m}_k.$
Note that the triple $(\mu^{(1)},$ $y^{(1)},$ $t^{(1)}$) is a two-cocycle. By our assumption, 
$$(\mu^{(1)},y^{(1)},t^{(1)})=f(\lambda^{(1)},x^{(1)},s^{(1)})+(d+\partial)(u^{(1)},r^{(1)})$$
for some cocycle $(\lambda^{(1)},x^{(1)},s^{(1)})$ and some cochain $(u^{(1)},r^{(1)}).$ Apply the gauge transformation $H={\operatorname{exp}}(u^{(1)}),$ $b={\operatorname{exp}}(r^{(1)})$ to $(\mu,$ $G,$ $c).$ We may assume that $(\mu^{(1)},y^{(1)},t^{(1)})=f(\lambda^{(1)},x^{(1)},s^{(1)})$ where $(\lambda^{(1)},x^{(1)},s^{(1)}$ is a cocycle.

By induction, we can replace $(\mu,$ $G,$ $c)$ by an isomorphic descent datum and assume that, modulo ${\mathfrak m}^{n+1},$ it is equal to $f(\lambda, F, a)$ where $(\lambda, G, a)$ is a descent datum modulo ${\mathfrak m}^{n+1}.$ By Corollary \ref{cor:non-abelian d2=0, II}, the cochain $(R^{(n+2)},$ $ Z^{(n+2)},$ $ T^{(n+2)},$ $ -\Phi^{(n+2)})$ is a cocycle. It is a coboundary, because its image under $f$ is (since $f(\lambda,F,a)$ is a descent datum), and $f$ is a quasi-isomorphism. Therefore, one can modify $(\lambda,F,a)$ in the component ${\mathfrak m}_{n+1},$ so that it will become a descent datum modulo ${\mathfrak m}^{n+2}.$ Furthermore, 
$$(d+\partial)(\mu^{(n+1)},y^{(n+1)},t^{(n+1)})-f(\lambda^{(n+1)}, x^{(n+1)},s^{(n+1)})=0,$$
therefore
$$(\mu^{(n+1)},y^{(n+1)},t^{(n+1)})-f(\lambda^{(n+1)}, x^{(n+1)},s^{(n+1)})=$$
$$(d+\partial)(u^{(n+1)}, r^{(n+1)})+f({\lambda'}^{(n+1)},{ x'}^{(n+1)},{s'}^{(n+1)})$$
where $({\lambda'}^{(n+1)},$ $ {x'}^{(n+1)},$ ${s'}^{(n+1)})$ is a cocycle. Replace $(\lambda^{(n+1)},$ $ x^{(n+1)},$ $s^{(n+1)}))$ by $(\lambda^{(n+1)}+{\lambda'}^{(n+1)},$ $ x^{(n+1)}+{x'}^{(n+1)},$ $s^{(n+1)}+{s'}^{(n+1)}),$ then apply the gauge transformation $H={\operatorname{exp}}(u^{(n+1)}),$ $b={\operatorname{exp}}(r^{(n+1)})$ to $(\mu,$ $G,$ $c.$ We get a new $(\lambda, G, a)$ which is a descent datum modulo ${\mathfrak m}^{n+2},$ and $f(\lambda, G, a)=(\mu, G, c)$  modulo ${\mathfrak m}^{n+2}.$ 

Now let us prove the injectivity in b). Let $(\lambda,$ $F,$ $a)$ and $(\lambda',$ $F',$ $a')$ be two descent data whose images under $f$ are isomorphic. Denote the isomorphism by $(H,$ $b).$ Let $F={\operatorname{exp}}(y),$ $F'={\operatorname{exp}}(y')$, $a={\operatorname{exp}}(s),$ $a'={\operatorname{exp}}(s'),$ $H={\operatorname{exp}}(x),$ $b={\operatorname{exp}}(r).$ We have
$$f(\lambda ^{(1)},y^{(1)}, s^{(1)})-f({\lambda'} ^{(1)},{y'}^{(1)},{ s'}^{(1)})=(d+\partial)(u^{(1)}, r^{(1)});$$
therefore, since $f$ is a quasi-isomorphism, the cocycle
$$(\lambda ^{(1)},y^{(1)}, s^{(1)})-({\lambda'} ^{(1)},{y'}^{(1)},{ s'}^{(1)})$$is a coboundary. After replacing the datum $(\lambda',$ $F',$ $a')$ by a datum which is isomorphic to it and identical to it modulo ${\mathfrak m}^2,$ we may assume that $f(\lambda,$ $F,$ $a)=f(\lambda',$ $F',$ $a')$ modulo ${\mathfrak m}^2.$ By induction, we may assume that $(\lambda,$ $F,$ $a)$ and $(\lambda',$ $F',$ $a')$ coincide modulo ${\mathfrak m}^{n+1}$ and that their images are isomorphic, the isomorphism being equal to identity modulo ${\mathfrak m}^{n}$.  Apply Corollary \ref{cor:non-abelian d2=0, III} to study the failure of $(H=1,$ $b=1)$ to be an isomorphism between $(\lambda,$ $F,$ $a)$ and $(\lambda',$ $F',$ $a')$. The corresponding cocycle is a coboundary because its image under $f$ is. Therefore we can act upon $(H=1,$ $b=1)$ by a two-morphism and obtain a new $(H,b)$ which is an isomorphism between  $(\lambda,$ $F,$ $a)$ and $(\lambda',$ $F',$ $a')$ modulo ${\mathfrak m}^{n+1}$. Now one can assume that $(\lambda,$ $F,$ $a)$ and $(\lambda',$ $F',$ $a')$ coincide modulo ${\mathfrak m}^{n+1}$ and that their images are isomorphic, the isomorphism being equal to identity modulo ${\mathfrak m}^{n+1}.$ We have $$f(\lambda ^{(n+1)},y^{(n+1)}, s^{(n+1)})-f({\lambda'} ^{(n+1)},{y'}^{(n+1)},{ s'}^{(n+1)})=(d+\partial)(u^{(n+1)}, r^{(n+1)});$$
since $f$ is a quasi-isomorphism, the cocycle
$(\lambda ^{(n+1)},y^{(n+1)}, s^{(n+1)})-({\lambda'} ^{(n+1)},{y'}^{(n+1)},{ s'}^{(n+1)})$ is a coboundary. After replacing the datum $(\lambda',$ $F',$ $a')$ by a datum which is isomorphic to it and identical to it modulo ${\mathfrak m}^{n+2},$ we may assume that $f(\lambda,$ $F,$ $a)=f(\lambda',$ $F',$ $a')$ modulo ${\mathfrak m}^{n+2}.$ 

This proves the statement b). The proofs of c) and d) are very similar, and we leave them to the reader.
\subsection{Totalization of cosimplicial DGLAs} \label{ss:Totalization of cosimplicial DGLAs} Here we recall how one can construct a DGLA from a cosimplicial DGLA by the procedure of totalization. We then prove that isomorphism classes of descent data for a cosimplicial DGLA are in one-to one correspondence with isomorphism classes of Maurer-Cartan elements of its totalization. This is a two-groupoid version of a theorem of Hinich \cite{H1}.

Define for $p\geq 0$
$${\mathbb Q}[\Delta ^p]={\mathbb Q}[t_{0}, \ldots, t_{p}]/(t_{0}+ \ldots + t_{p}-1)$$
and
$$\Omega ^{\bullet} [\Delta  ^p]={\mathbb Q}[t_{0}, \ldots, t_{p}]\{dt_{0}, \ldots, dt_{p}\}/(t_{0}+ \ldots + t_{p}-1,dt_{0}+ \ldots + dt_{p} )$$
The collection $\{\Omega ^{\bullet} [\Delta  ^p]\}, p\geq 0,$ is a simplicial DGA. 

Let $\cM$ be the category whose objects are morphisms $f:[p]\to [q]$ in $\Delta$ and a morphism from $f:[p]\to [q]$ to $f':[p']\to [q']$ is a pair $a: [p']\to [p]$, $b:[q]\to [q']$ such that $f'=bfa.$ Given $(a,b):f\to f'$ and $(a',b'): f'\to f'',$ define their composition to be $(a'a,bb').$ 

Given a cosimplicial DGLA $\cL,$ we can construct a functor from $\cM$ to the category of vector spaces by assigning to the object $f:[p]\to [q]$the space $\Omega ^{\bullet} [\Delta  ^p]\otimes \cL^q.$ Set 
$$\Tot (\cL)=\limdir _{\cM} \Omega ^{\bullet} [\Delta  ^p]\otimes \cL^q$$
This is a DGLA (with the differential being induced by $d_{\operatorname{DR}}.$
\begin{proposition}\label{prop:Hinich}
a). There is a bijection between the set of isomorphism classes of descent data of the Deligne two-groupoid of $\cL$ and the set of Maurer-Cartan elements of $\Tot(\cL).$ 

b). For a descent datum $\cA$ of $\cL$, denote by $\lambda(\cA)$ a Maurer-Cartan element from the isomorphism class given by a). Then there is a bijection
$$\frac{{\operatorname{Iso}}(\cA, \cA')}{2-{\operatorname{Iso}}}\isomoto \frac{{\operatorname{Iso}}(\lambda(\cA), \lambda(\cA'))}{2-{\operatorname{Iso}}}.$$
c). For two isomorphisms $\phi, \psi:\cA\to \cA'$, denote their images under the above bijection by $G(\phi), G(\psi).$ Then $f$ induces a bijection
$$2-{\operatorname{Iso}}(\phi, \psi)\isomoto 2-{\operatorname{Iso}}(G(\phi), G(\psi))$$
\end{proposition}
{\bf Proof.}  Recall that for every small category $\cM$ and for every functor $C:\cM\to {\operatorname{Vect}}_k$ one can define a cosimplicial space 
$$({\bf R}\liminv _{\cM}C)^n=\prod _{f_0\stackrel{\alpha_1}{\rightarrow} f_1\stackrel{\alpha_2}{\rightarrow}\ldots  \stackrel{\alpha_n}{\rightarrow}f_n} C(f_n)$$
with the standard maps $d_i$ and $s_i.$ The product is taken over all composable chains of morphisms in $\cM.$ If $C$ is a functor from $\cM$ to the category of $DGLAs$ then ${\bf R}\liminv _{\cM}C$ is a cosimplicial DGLA. 
 
 Consider the cosimplicial DGLA
${\bf R}\liminv _{\cM}\Omega ^{\bullet} [\Delta  ^p]\otimes \cL^q,$ together with the constant cosimplicial DGLA $\Tot(\cL)$  and the cosimplicial DGLA ${\bf R}\liminv _{\Delta}(\cL)$. The second and the third DGLAs embed into the first, and these embeddings are quasi-isomorphisms with respect to the differentials $d+\partial.$ By Proposition \ref{prop:quis of cosimplicial dglas}, our statement is true if we replace the cosimplicial DGLA $\cL$ by ${\bf R}\liminv _{\Delta}(\cL)$. But these two cosimplicial DGLA are quasi-isomorphic, whence the statement.
\subsection{The Hochschild complex} \begin{definition} \label{ex:Hochschild complex} For any associative algebra $A$, let $\cL ^H (A)$ be the Hochschild cochain complex equipped with the Gerstenhaber bracket \cite{Ge}. The standard Hochschild differential is denoted by $\delta$. For a sheaf of algebras $\cA,$ let $\cL ^H(\cA)$ denote the sheafification of the presheaf of DGLA $U\mapsto \cL ^H(A(U)).$ For the sheaf of algebras $C^{\infty}_M$ on a smooth manifold, resp. ${\cal O}_M$ on a complex analytic manifold, let ${\cL} ^H_M$ be the sheaf of Hochschild cochains $D(f_1, \ldots, f_n)$ which are given by multi-differential, resp. holomorphic multi-differential, expressions in $f_1, \ldots, f_n$.\end{definition}

One gets directly from the definitions the following

\begin{lemma} \label{lemma:defs of triv gerbe through DGLA} 
The set of isomorphism classes of deformations over ${\mathfrak{a}}$ of a sheaf of $k$-algebras $\cA$ as a stack is in one-to-one correspondence with the set of isomorphism classes of descent data of the Deligne two-groupoid of $\cL ^H(\cA)\otimes{\mathfrak{m}}$.  Similarly, the set of isomorphism classes of deformations of the trivial gerbe on $M$ is in one-to-one correspondence with the set of isomorphism classes of descent data of the Deligne two-groupoid of $\cL ^H _M\otimes{\mathfrak{m}}$.
\end{lemma}

\subsection{Hochschild cochains at the jet level} \label{ss:Hochschild cochains at the jet level} For a manifold $M$, let $J$, or $J_M$, be the bundle of jets of smooth, resp. holomorphic, functions on $M$. By $\nabla _{{\operatorname {can}}}$ we denote the canonical flat connection on the bundle $J$. Let $C^{\bullet}(J,J)$ be the bundle of Hochschild cochain complexes of $J$. More precisely, the fibre of this bundle is the complex of jets of multi-differential multi-linear expressions $D(f_1, \ldots, f_n)$. We denote by $\delta$ the standard Hochschild differential. 

\begin{proposition} \label{prop:defs of triv gerbe through DGLA of jets} The set of isomorphism classes of deformations of the trivial gerbe on $M$ is in one-to-one correspondence with the set of isomorphism classes of Maurer-Cartan elements of the DGLA $\cL ^{H,J} (M)\otimes{\mathfrak m}$ where
$$\cL ^{H,J} (M)= A ^{\bullet} (M,C^{\bullet +1} (J,J)) $$ with the differential $\nc + \delta$.
Here by $A^{\bullet}$ we mean $C^{\infty}$ forms with coefficients in a bundle.
\end{proposition}
{\bf Proof.} We have an embedding of sheaves of DGLA:
$$\cL ^H_M \to A ^{\bullet} _M (C^{\bullet +1} (J,J)) $$
which is a quasi-isomorphism, and the sheaf on the right hand side has zero cohomology in positive degrees. The proposition follows from Proposition \ref{prop:quis of cosimplicial dglas}.

\section{Deformation quantization of the trivial gerbe on a symplectic manifold} \label{s:Deformations of the trivial gerbe on a symplectic manifold} 
 \subsection{Deformation quantization of gerbes}
\begin{definition}\label{dfn:deformation quantization} A deformation quantization of a gerbe $\cA ^{(0)}$ on a manifold $M$ is a collection of deformations $\cA^{(N)}$ over ${\mathfrak{a}}={\mathbb{C}}[\hbar]/(\hbar^{N+1}),$ $N\geq 0$ (cf. Definition \ref{dfn:deformation1}), such that $\cA^{(N)}/\hbar^{N}=\cA^{(N-1)}.$ An isomorphism of two deformation quantizations is a collection of isomorphisms of deformations $\varphi _N:\cA^{(N)}\to {\cA'}^{(N)}$ such that $\varphi _N \congr \varphi _{N-1} \;{\operatorname{mod}}\hbar ^N.$
\end{definition}
Given a deformation quantization of a gerbe, one can define a stack of ${\mathbb{C}}[\hbar]$-algebras $\cA=\liminv \cA^{(N)}.$ Usually we will not distinguish between the deformation quantization and this stack.
\subsection{} Let $(M, \omega)$ be a symplectic manifold ($C^{\infty}$ or complex analytic with a holomorphic symplectic form). In this section, we extend Fedosov's methods from \cite{Fe} to deformations of the trivial gerbe. We say that a deformation quantization of the trivial gerbe on $M$ corresponds to $\omega$ if, on every $U_k$, $f*g-g*f = {\sqrt{-1}}\hbar \{f,g\} + o(\hbar)$ where $\{\,,\,\}$ is the Poisson bracket corresponding to $\omega$.

Let us observe that the group $H^2 (M, \hbar \C[[\hbar]])$ acts on the set of equivalence classes of deformations of any stack: a class $\gamma$ acts by multiplying $c_{ijk}$ by ${\operatorname {exp}}\,\gamma _{ijk}$ where $\gamma _{ijk}$ is a cocycle representing $\gamma$.
\begin{thm} \label{thm:classification of deformations of the trivial gerbe, symplectic case} Denote by $\operatorname{Def}(M,\omega)$ the set of isomorphism classes of deformation quantizations of the trivial gerbe on $M$ compatible with the symplectic structure $\omega$. The action of $H^2 (M, \hbar \C[[\hbar]])$ on $\operatorname{Def}(M,\omega)$ is free. The space of orbits of this action is in one-to-one correspondence with an affine space modelled on the vector space $ H^2 (M, \C)$ (in the $C^{\infty}$ case) or $H^1(M,{\cO _M}/{\Bbb C})$ (in the complex case).
\end{thm}
{\bf Proof.} As in \cite{Fe1}, we will reduce the proof to a classification problem for certain connections in an infinite-dimensional bundle of algebras. 

Let us observe that the Proposition \ref{prop:defs of triv gerbe through DGLA of jets} is true if we replace deformations over Artinian rings by deformation quantizations. Indeed, the proof of Proposition \ref{prop:quis of cosimplicial dglas} works verbatim for the DGLAs that are needed for Proposition \ref{prop:defs of triv gerbe through DGLA of jets}, since one can start with a good cover, and all cohomological obstructions are zero already in the \v{C}ech complex of this cover; one has no need of refining the cover, and therefore one can carry out the induction procedure infinitely many times. Next, note that in Proposition \ref{prop:defs of triv gerbe through DGLA of jets} we can replace the bundle of algebras $J$ by the bundle of algebras 
$${\operatorname{gr}} J = \prod S^m(T^*_M).$$
 Indeed, a standard argument shows that they are isomorphic as $C^{\infty}$ bundles of algebras. 

Under this isomorphism, the canonical connection $\nc$ becomes a connection $\nabla _0$ on ${\operatorname{gr}} J$. We are reduced to classifying up to isomorphism those Maurer-Cartan elements of $(A^{\bullet}(M, C^{\bullet +1}({\operatorname{gr}} J, {\operatorname{gr}} J)), \nabla _0 +\delta)$ whose component in $A^0(M, C^2)$ is equal to $\frac{1}{2}{\sqrt{-1}}\hbar \{f,g\}$ modulo $\hbar$. In other words,these components must be, pointwise, deformation quantizations of $\prod S^m(T^*_M)$ corresponding to the symplectic structure. But all such deformations are isomorphic to the standard Weyl deformation from the definition below:
\begin{definition} \label{dfn:Weyl algebra}
The Weyl algebra of $T^*_M$ is the bundle of algebras 
$$W = {\operatorname{gr}} J[[\hbar]] = \prod S^m(T^*_M)[[\hbar]]$$
with the standard Weyl product $*$.
\end{definition}
 Moreover, a smooth field of such deformations on $M$ admits a smooth gauge transformation making it the standard Weyl deformation. Therefore, we have to classify up to isomorphism those Maurer-Cartan elements of $A^{\bullet}(M, C^{\bullet +1}({\operatorname{gr}} J, {\operatorname{gr}} J))$ whose component in the subspace $A^0(M, C^2)$ is equal to $f*g - fg$. Here $*$ is the product in the standard Weyl deformation.

\begin{proposition}\label{prop:Fedosov for stacks}
Deformations of the trivial gerbe on $M$ compatible with a symplectic structure $\omega$ are classified up to isomorphism 
by pairs $(A,c)$ where 
\begin{equation}\label{eq:Fedosov 1}
A\in \hbar A^1 (M, {\operatorname{hom}}({\operatorname{gr}} J, {\operatorname{gr}} J))[[\hbar]];
\end{equation}
\begin{equation}\label{eq:Fedosov 2}
c\in \hbar A^2 (M, {\operatorname{gr}} J)[[\hbar]],
\end{equation}
such that, if 
$$\nabla = \nabla _0 + A,$$
then
\begin{equation}\label{eq:Fedosov 3}
\nabla (f*g) = \nabla(f) *g + f*\nabla(g);
\end{equation}
\begin{equation}\label{eq:Fedosov 4}
\nabla ^2 = {\operatorname{ad}}(c);\;\nabla (c) = 0
\end{equation}
Two pairs $(A,c)$ and $(A',c')$ are equivalent if one is obtained from the other by a composition of transformations of the following two types.
a)
\begin{equation}\label{eq:Fedosov 5}
 (A,\,c) \mapsto ({\operatorname{exp}}({\operatorname{ad}}(X))(A),\, {\operatorname{exp}}({\operatorname{ad}}(X))(c))
\end{equation}
where $X \in \hbar{\operatorname{Der}}(W);$

b)
\begin{equation}\label{eq:Fedosov 6}
 (A,\,c) \mapsto (A+B,\,c+{\nabla}B +{\frac{1}{2}}[B,\,B])
\end{equation}
where $B \in \hbar W$.
\end{proposition}

It is straightforward that the set of Maurer-Cartan elements discussed above, up to isomorphism, is in one-to-one correspondence with the set of pairs $(A,\,c)$ up to equivalence. Indeed, given $(A,c)$, the Maurer-Cartan element is constructed as follows: the component in $A^0(M, C^2)$ is the difference between the Weyl product and the commutative product; the component in $A^1(M, C^1)$ is $\nabla-\nabla _0$, and  the component in $A^2(M, C^0)$ is $c$. It remains to show that the pairs $(A,\,c)$ are classified as in Theorem \ref{thm:classification of deformations of the trivial gerbe, symplectic case}.

Let us start with notation. Let 
$${\widetilde{\g}}^0={\operatorname{ gr}}{J}$$
be the bundle of Lie algebras of formal power series with the standard Poisson bracket. Let $\g ^0 = {\operatorname{ gr}}{J}/{\Bbb C}$ 
be the quotient bundle of Lie algebras. In other words, the fibre of $\g ^0$ is the Lie algebra of formal Hamiltonian vector fields on the tangent space. Also, put
$${\widetilde{\g}}={\frac{1}{\hbar}}W$$
with the bracket $a*b-b*a$ where $*$ is the Weyl product, and
$$\g = {\widetilde{\g}}/{\frac{1}{\hbar}}{\Bbb C}[[\hbar]]$$
This is the Lie algebra of continuous derivations of the Weyl algebra. It maps surjectively to $\g ^0$ via ${\frac{1}{\hbar}}(f_0 + \hbar f_1 + \cdots)\mapsto f_0$. Put $|a|=m$ for $a\in S^m(T^*_M)$ and $|\hbar|=2$. This defines the degree of any monomial in $S^m(T^*_M)[\hbar].$ By ${\widetilde{\g}}^{0}_m$ we denote the subspace $S^{m+2}(T^*_M)$, and by ${{\widetilde{\g}}}_{m}$ the set of ${\frac{1}{\hbar}}f$ where $f$ is a polynomial from $S^{m+2}(T^*_M)[\hbar].$ Then 
$$[{\widetilde{\g}}^0_{m},\,{\widetilde{\g}}^{0}_{r}]\subset {\widetilde{\g}}^{0}_{m+r};\;\;[{\widetilde{\g}}_{m},\,{\widetilde{\g}}_{r}]\subset {\widetilde{\g}}_{m+r};$$
$${\widetilde{\g}}^{0}=\prod_{m\geq -2} {\widetilde{\g}}^0_{m};\;{\widetilde{\g}}=\prod_{m\geq -2} {\widetilde{\g}}_{m}$$
One defines $\g ^0_{m}$ and $\g _m$ accordingly. We have
$${{\g}}^{0}=\prod_{m\geq -1} {{\g}}^{0}_{m};\;{{\g}}=\prod_{m\geq -1} {{\g}}_{m}$$
 In particular, the bundle ${\widetilde{g}}^0_{-1}=\g^{0}_{-1}={\widetilde{g}}_{-1}=\g_{-1}$ is the cotangent bundle $T^*_M$. The symplectic form identifies this bundle with $T_M$.
\begin{definition}\label{df:A_-1}
By $A_{-1} $ we denote the canonical form ${\operatorname{id}}\in A^1(M,T_M)$ which we view as a form with values in ${\widetilde{g}}^{0}_{-1},$ etc. under the identifications above.
\end{definition}
The form $A_{-1}$ is smooth in the $C^{\infty}$ case and holomorphic in the complex case.

The connection $\nabla _0$ can be expressed as
\begin{equation} \label{eq:nabla_0}
\nabla _0=A_{-1}+\nabla _{0,0}+\sum_{k=1}^{\infty} A_k=\nabla _{0,0}
\end{equation}
where $\nabla _{0,0}$ is an ${\frak{sp}}_n$-valued connection in the tangent bundle $T_M$ and $A_k \in A^1(M, \g^{0}_{k})$. 
Define
$$A^{(-1)}=\sum_{k=1}^{\infty} A_k$$
(Here $n=\frac{1}{2}{\operatorname{dim}}(M)$).
The form $A_{-1}$ is in fact the canonical form from the above definition. In the case of a complex manifold, locally $\nabla _{0,0}=\partial +{\overline {\partial}}+A_{0,0}$ where $A_{0,0}$ is a $(1,0)$-form with values in ${\frak{sp}}_n$. The form $A^{(-1)}$ can be viewed as a ${\widetilde{\g}}^0$-valued one-form:
\begin{equation} \label{eq:nabla_0, I}
A^{(-1)}\in A^1(M, {\widetilde{\g}}^0)
\end{equation}
Let us look for $\nabla$ of the form
\begin{equation}\label{eq:nabla in general}
\nabla=\no + \sum_{m=0}^{\infty}({\sqrt{-1}}\hbar)^m A^{(m)}
\end{equation}
where $A^{(m)} \in A^1(M, \g ^0)$. The condition $\nabla ^2 = o(\hbar)$ is equivalent to 
\begin{equation}\label{eq:flatness mod h}
\no A^{(0)}+{\frac{1}{2}}[A^{(-1)},A^{(-1)}]_{2}=0 
\end{equation}
Here we use the notation
$$a*b-b*a=\sum_{m=1}^{\infty} ({\sqrt{-1}}\hbar)^m [a,b]_{m}$$
(in particular, $[\,,\,]_0$ is the Poisson bracket); we then extend the brackets $[a,b]_{m}$ to forms with values in the Weyl algebra. Since $[\nc,[\nc,\nc]]=0$ and $[\no,\no]=0,$ we conclude that 
$$\no [A^{(-1)}, A^{(-1)}]_2 =0$$
in $A^2 (M, {\widetilde{\g}}^0)$. Moreover, observe that the left hand side lies in fact in  $A^2 (M, \prod _{m\geq 0}{\widetilde{\g}}^0_{m}).$ 
\begin{lemma}\label{lemma:acyclicity of A-1}
If $c\in A^p(M, {\widetilde{\g}}^0_{m})$, $m\geq -1,$ satisfies $[A_{-1},c]=0,$ then $c=[A_{-1},c']$ for  $c'\in A^{p-1}(M, {\widetilde{\g}}^0_{m+1}).$
\end{lemma}
{\bf Proof.} Indeed, the complex $A^{\bullet}(M, {\widetilde \g}^0)$ with the differential $[A_{-1},\;]$ is isomorphic to the complex of smooth sections of, resp, $A^{0,{\bullet}}$ forms with coefficients in, the bundle of complexes $S[[T^*_M]]\otimes \wedge (T^*_M)$ with the standard De Rham differential.

We now know that pairs $(\nabla, c)$ exist. The theorem is implied by the following lemma (we use the notation of \eqref{eq:Fedosov 1}-\eqref{eq:Fedosov 6}).
\begin{lemma} \label{lemma:classification of connections} 
1) For any two connections $\nabla$ and $\nabla ',$ $A^{(0)}-{A'}^{(0)}$ is a cocycle in $A^1(M, J/{\C});$ a pair $(\nabla, c)$ is equivalent to a pair $(\nabla ', c')$ for some $c'$ by some transformation $(X, B)$ if and only if $A^{(0)}-{A'}^{(0)}$ is a coboundary;

2) for any two pairs $(\nabla, c)$ and $(\nabla, c')$ with the same $\nabla,$ $c-c'$ is a closed form  in $A^2 (M, \hbar \C[[\hbar]]);$ two such pairs are equivalent if and only if $c-c'$ is exact.
\end{lemma}
{\bf Proof.} 1) The first statement of 1) follows from \eqref{eq:flatness mod h}. To prove the second, note that 
$$\nabla ' = \ea (X) (\nabla) + {\operatorname ad}(B),$$
$$B \in A^1 (M, \hbar {\widetilde {\g}})$$
with
$$X=\sum_{m=0}^{\infty}({\sqrt {-1}}\hbar)^m X^{(m)}$$
and $X^{(m)}\in A^0(M, \g ^0),$ is possible if and only if 
$$\no X^{(0)} + A^{({0})}-{A'}^{({0})}=0.$$
2) The first statement of 2) follows from \eqref{eq:Fedosov 4}. To prove the second, consider a lifting of $\nabla$ to a ${\widetilde{g}}$-valued connection $\tn.$ We have
$$c=\tn ^2 + \theta$$
where $\theta \in A^2 (M, \hbar \C[[\hbar]])$. One has 
$$\nabla = \ea (X)(\nabla)+B$$
if and only if the following two equalities hold:
$$\tn = \ea (X)(\tn)+B+\alpha$$
for some $\alpha \in A^1 (M, \C[[\hbar]]);$
$$c' = \ea (X)(c)+\ea (X)(B)+{\frac{1}{2}}[B,B].$$
But in this case
$$c'=\ea (X) (\tn ^2 + \theta)+[\ea (X)(\tn), \tn-\ea (X)(\tn)-\alpha]+$$
$${\frac{1}{2}}[\tn-\ea (X)(\tn),\tn-\ea (X)(\tn)]=$$
$${\frac{1}{2}}[\ea (X)(\tn),\ea (X)(\tn)]+$$
$$\theta + [\ea (X) \tn, \tn]- {\frac{1}{2}}[\ea (X)(\tn),\ea (X)(\tn)] -d\alpha +{\frac{1}{2}}[\tn,\tn]-$$
$$[\tn, \ea (X)(\tn)]+
{\frac{1}{2}}[\ea (X)(\tn),\ea (X)(\tn)]=\tn ^2 +\theta -d\alpha$$
$$ =c-d\alpha$$
This proves the theorem.

\section{The characteristic class of a deformation and the Rozansky-Witten class} \subsection{The characteristic class} Given a deformation of the trivial gerbe on a symplectic manifold $(M,\omega)$, one defines its characteristic class 
$$\theta= \frac{1}{{\sqrt{-1}}\hbar}\omega+\sum_{k=0}^{\infty}({{\sqrt{-1}}\hbar})^k\theta_k\in \frac{1}{{\sqrt{-1}}\hbar}\omega+H^2(M)[[\hbar]]$$
as follows. Represent the deformation by a pair $(\nabla, c)$ as in Proposition \ref{prop:Fedosov for stacks}. Choose a lifting $\tn$ of $\nabla$ to a ${\mathfrak{g}}$-valued connection; define
$$\theta =\tn ^2-c.$$
It is easy to see that:
i) $\theta\in A^2(M, {\frac{1}{\hbar}}{\mathbb C}[[\hbar]]);$

ii) $d\theta=0,$ and the cohomology class of $\theta$ is invariant under the equivalence and independent of the lifting. 

The above construction generalizes Fedosov's Weyl curvature. It is easy to see that the class of$ \theta_0$ coincides with the image of the class from Theorem \ref{thm:classification of deformations of the trivial gerbe, symplectic case} under the morphism
$\partial: H^1(M, \cO_M/{\mathbb C})\to H^2(M, {\mathbb C}).$ In particular, if this map is not injective, there may be non-isomorphic deformations with the same class $\theta$.
\subsection{Deformation quantization of the sheaf of functions}
Here we recall a theorem from \cite{NT} (cf. \cite{BK} for the algebraic case).

Let $(M, \omega)$ be either a symplectic $C^{\infty}$ manifold or a complex manifold with a holomorphic symplectic structure. By $\cO _M$ we denote the sheaf of smooth, resp. holomorphic, functions.

In what follows we will study deformation quantization of $\cO_M$ {\em as a sheaf}. In the language adopted in this article, these are deformation quantizations of the trivial gerbe such that $c_{ijk}=1$. An isomorphism is by definition an isomorphism of deformation quantizations such that $b_{ij}=1$.
\begin{thm} \label{thm:classification of sheaves} 
Assume that the maps $H^i (M, \C) \to H^i (M, \cO _M)$ are onto for $i=1,\,2$. 
Set 
$$H^2_F (M,\C)={\operatorname {ker}}(H^2(M, \C) \to H^2(M, \cO_M)).$$
Choose a splitting
$$H^2(M, \C) = H^2(M, \cO_M) \oplus H^2_F (M,\C).$$
The set of isomorphism classes of deformation quantizations of $\cO _M$ as a sheaf which are compatible with $\omega$ is in one-to-one correspondence with a subset of the affine space 
$${\frac{1}{{\sqrt{-1}}\hbar}}\omega + H^2 (M, \C)[[\hbar]]$$ 
whose projection to
$${\frac{1}{{\sqrt{-1}}\hbar}}\omega + H^2 _F (M, \C)[[\hbar]]$$
is a bijection.
\end{thm}

\subsection{The first Rozansky-Witten class} \label{s:the RW class}
 We have seen in the previous section that, under the assumptions of Theorem \ref{thm:classification of sheaves}, deformations of the sheaf of algebras $\cO _M$ are classified by cohomology classes $\theta$ as in \eqref{eq:teta} where $\theta _{-1}=\frac{1}{{\sqrt {-1}}\hbar }\omega;$ the (non-natural) projection of the set of all possible classes $\theta$ to $\frac{1}{{\sqrt {-1}}\hbar }\omega + H^2_F(M,\C[[\hbar]])$ is a bijection. More precisely, the (natural) projection of $\theta _{n+1}$ to $H^2 (M, \cO_M)$ is a nonlinear function in $\theta _i ,\;0\leq i \leq n.$ We are going to describe this function for the case $n=0.$

Let $M$ be a complex manifold with a holomorphic symplectic structure $\omega$. We start by describing two ways of constructing cohomology classes in $H^2 (M, \cO _M) .$ The first one was invented by Rozansky and Witten, cf. \cite{RW}, \cite{Kap}, \cite{K2}. Let $\nabla _{0,0}$ be a torsion-free connection in the tangent bundle which is locally of the form $d+A_0$ for $A_0 \in A^{1,0}(M, {\mathfrak{sp}}).$
 Let $R={\overline{\partial}}A_0$ be the $(1,1)$ component of the curvature of $\nabla _{0,0}$. We can view $R$ as a $(1,1)$ form with coefficients in $S^2 (T^*_M)$. Let $z^i$ be holomorphic coordinates on $M.$ By ${\widehat{z}}^i$ we denote the corresponding basis of $T^*_M.$ We write
\begin{equation}\label{eq:R}
R=\sum R_{abi{\overline{j}}}{\widehat{z}}^a{\widehat{z}}^bdz^id{\overline{z}}^j
\end{equation}
Put 
\begin{equation}\label{eq:RW}
{\operatorname{RW}}_{\Gamma _0}(M,\omega)=\sum R_{abi{\overline{j}}}  R_{cdk{\overline{l}}}\omega^{ac}\omega^{bd}\omega^{ik}d{\overline{z}}^jd{\overline{z}}^l
\end{equation}
Here $\Gamma _0$ refers to the graph with two vertices and three edges connecting them. In fact a similar form ${\operatorname{RW}}_{\Gamma }(M,\omega)$ can be defined for any finite graph $\Gamma$ for which every vertex is adjacent to three edges; the cohomology class of this form is independent of the connection \cite{RW}.

The other way of obtaining $(0,2)$ classes is as follows. For $\alpha = \sum \alpha _{i{\overline{j}}}dz^i d{\overline{z}}^j$ and $\beta = \sum \beta _{i{\overline{j}}}dz^i d{\overline{z}}^j,$ put
\begin{equation}\label{eq:def of pairing}
\omega(\alpha, \beta)=\sum \alpha _{i{\overline{j}}} \beta _{k{\overline{l}}}\omega_{ik}d{\overline{z}}^jd{\overline{z}}^l
\end{equation}
It is straightforward that the above operation defines a symmetric pairing 
$$\omega: H^{1,1}(M)\otimes H^{1,1}(M) \to H^{0,2}(M).$$
Combined with the projection $H^2_F(M)\to H^{1,1}(M),$ this gives a symmetric pairing
$$\omega: H^2_F(M)\otimes H^2_F(M) \to H^2(M,\cO_M).$$
\begin{thm} \label{thm:RW}
Under the assumptions of Theorem \ref{thm:classification of sheaves}, let a deformation of the sheaf of algebras $\cO_M$ correspond to a cohomology class 
$$\theta=\sum ({\sqrt{-1}}\hbar)^m \theta _m, \; \theta _m \in H^2 (M).$$
 Then the projection of the class of $\theta _1$ to $H^2 (M,\cO _M)$ is equal to
$${\operatorname{RW}}_{\Gamma _0}(M,\omega)+\omega(\theta _0,\theta _0)$$
\end{thm}
{\bf Proof.} 
 First, observe that Lemma \ref{lemma:defs of triv gerbe through DGLA} and Proposition \ref{prop:defs of triv gerbe through DGLA of jets} have their analogs for deformations of the structure sheaf as a sheaf of algebras. The only difference is that the Hochschild complex $C^{\bullet +1}$ is replaced everywhere by $C^{\bullet +1}, {\bullet \geq 0}.$ Similarly to \eqref{eq:Fedosov 1}-\eqref{eq:Fedosov 6}, one has
\begin{lemma} \label{lemma:Fedosov for sheaves}  Deformations of the sheaf of algebras ${\cal{O}}_M$ which are compatible with a symplectic structure $\omega$ are classified by forms 
$A\in \hbar A^1 (M, {\operatorname{hom}}({\operatorname{gr}} J, {\operatorname{gr}} J))[[\hbar]]$
such that, if 
$$\nabla = \nabla _0 + A,$$
then
\begin{equation}\label{eq:Fedosov 33}
\nabla (f*g) = \nabla(f) *g + f*\nabla(g);
\end{equation}
and $\nabla ^2 = 0$. Two such forms are equivalent if, for $X\in  A^0 (M, \hbar {\operatorname{Der}}(W),$ 
$$\nabla ' =\ea (X)\nabla$$
\end{lemma}
The proof is identical to the proof of Lemma \ref{prop:Fedosov for stacks}. 

Let us now classify pairs $(\nabla, c)$. 

We start by constructing a flat connection $\nabla$. We use a standard proof from the homological perturbation theory. One has to solve recursively
\begin{equation} \label{eq:curvature}
R_n +\no A^{(n+1)}=0
\end{equation}
where
$$R_n=\frac{1}{2}\sum_{i,j\geq 0\;;i+j+m=n+1} [A^{(i)},A^{(j)}]_m $$
At every stage $\no R_n = 0$; the class of $R_n$ is in the image of the map 
$$H^2(M, \cO _M)\to H^2(M, \cO _M /\C)$$ which is zero under our assumptions.

We have shown that flat connections $\nabla$ exist. For any such connection we can consider its lifting to a ${\widetilde {\g}}$-valued connection $\tn .$ Put
\begin{equation} \label{eq:theta}
\tn ^2 = \theta =  \sum _{m=-1}^{\infty}({\sqrt{-1}}\hbar)^m \theta _m\in A^2 (M, {\frac {1}{\hbar}}\C[[\hbar]])
\end{equation}
Let us try to determine all possible values of $\theta.$ 
\begin{lemma}\label{lemma:nablas via tetas}
Under the assumptions of Theorem \ref{thm:classification of sheaves}, the map $\nabla\mapsto \theta$ establishes a one-to-one correspondence between the set of equivalence classes of connections $\nabla$ and a subset of the affine space 
$${\frac{1}{{\sqrt{-1}}\hbar}}\omega + H^2 (M, \C)[[\hbar]]$$ 
whose projection to
$${\frac{1}{{\sqrt{-1}}\hbar}}\omega + H^2 _F (M, \C)[[\hbar]]$$
is a bijection.
\end{lemma}First of all, $\theta _{-1}=\frac{1}{\sqrt{-1}\hbar}\omega.$ There exists $\tn$ with $\theta _0 =0$ (see \eqref{eq:flatness mod h} and the argument after it). To obtain other possible $\theta _0$ we have to add to $\tn$ a form ${A'}^{(0)}-{A}^{(0)}$ whose image in $A^1(M, J/\C)$ is $\tn$-closed. Therefore, the cohomology class of a possible $\theta _0$ must be in the image of the map
$$H^1 (M, \cO _M/\C)\to H^2 (M, \C),$$
which is precisely $H^2_F(M, \C)$ under our assumptions.

Proceeding by induction, we see that, having constructed $\theta _i, $ $i\leq n,$ and $\tn _{(n)}$ such that 
\begin{equation}\label{eq:teta}
\tn _{(n)}^2  =  \sum _{m=-1}^{n}({\sqrt{-1}}\hbar)^m \theta _m +o(\hbar ^n),
\end{equation}
we can find $\theta _{n+1}$ and $\tn _{(n+1)}=\tn ^{(n)}+o(\hbar ^n)$ such that
$$\tn _{(n+1)}^2  =  \sum _{m=-1}^{n+1}({\sqrt{-1}}\hbar)^m \theta _m +o(\hbar ^{n+1}).$$
The cohomology class of such $\theta_{n+1}$ can be changed by adding any element of $H^2_F(M)$. 

Proceeding by induction, we see that we can construct unique $\tn$ with any given projection of $\theta$ to $H^2_F(M)[[\hbar]].$ Now observe that, if $\nabla' = \ea (X)\nabla,$ then $\tn ' = \ea (X)\tn + \alpha$ for $\alpha \in A^1 (M, \C[[\hbar]])$ and therefore $\theta ' = \ea(X) (\theta) + d\alpha.$ Therefore two connections with non-cohomologous curvatures are not equivalent. An inductive argument, similar to the ones above, shows that two connections with cohomologous curvatures are equivalent. Indeed, by adding an $\alpha$ we can arrange for $\theta '$ and $\theta$ to be equal. Then we find $X=\sum ({\sqrt {-1}}\hbar)^m X_m$ by induction. At each stage we will have an obstruction in the image of the map 
$$H^1(M,\cO _M)\to H^1(M,\cO _M/\C).$$
But this image is zero under our assumptions. 

\subsubsection{End of the proof of Theorem \ref{thm:RW} }Let us start by observing that one can define the projection
\begin{equation}\label{eq:Proj}
{\operatorname {Proj}}:( A^{\bullet, \bullet}(M, {\operatorname {gr}}\,J), \no)  \to (A^{0, \bullet}(M), {\overline{\partial}})
\end{equation}
as follows: if ${\cal I}$ is the DG ideal of the left hand side generated by $dz ^i$ and by the augmentation ideal of ${\operatorname {gr}}\,J$ then the right hand side is identified with the quotient of the left hand side by ${\cal I.}$ It is straightforward that ${\operatorname {Proj}}$ is a quasi-isomorphism.

Using the notation introduced in and after Definition \ref{df:A_-1}, we can write
\begin{equation}\label{eq:flatness mod h again}
\no A^{(0)}+{\frac{1}{2}}[A^{(-1)},A^{(-1)}]_{2}=\theta _0 
\end{equation}
and
\begin{equation}\label{eq:flatness mod h^2}
\no A^{(1)}+{\frac{1}{2}}[A^{(-1)},A^{(-1)}]_{3}+[A^{(-1)},A^{(0)}]_{2}+[A^{(-0)},A^{(0)}]_{1}=\theta _1.
\end{equation}
Observe that:

a)
${\operatorname {Proj}}[A^{(-1)},A^{(-1)}]_{2}={\operatorname {Proj}}[A^{(-1)},A^{(-0)}]_{2}=0;$

b) ${\operatorname {Proj}}[A^{(-1)},A^{(-1)}]_{3}$ depends only on the $(0,1)$ component of the form $A^{(-1)}_1;$

c) ${\operatorname {Proj}}[A^{(0)},A^{(0)}]_{1}$ depends only on the $(0,1)$ component of the form $A^{(0)}_{-1}.$

The connection $\no$ can be chosen in such a way that the form from b) is equal to
\begin{equation}\label{eq:A-1+101}
\sum R_{ijk{\overline{l}}}{\widehat{z}}^i{\widehat{z}}^j{\widehat{z}}^kd{\overline{z}}^l;
\end{equation}
therefore for this connection
 $${\frac{1}{2}}{\operatorname {Proj}}[A^{(-1)},A^{(-1)}]_{3}={\operatorname RW}_{\Gamma _0}(M,\omega).$$
Since $[A^{(-1)},A^{(-1)}]_{2}\in A^2 (M,{\widetilde{\g}}_{\geq 0}),$ we can choose $A^{(0)}\in A^1 (M,{\widetilde{\g}}_{\geq 1});$ we conclude, because of b) and c), that there exists $\tn$ with $\theta _0=0$ such that the projection of $\theta _1$ to $H^2 (M, \cO_M)$ is equal to ${\operatorname RW}_{\Gamma _0}(M,\omega).$

Now we can produce a connection with a given $\theta _0$ by adding to the above connection a form $A'-A;$ for this new connection, the form from c) may be chosen as 
$$\sum \alpha _{i{\overline j}}{\widehat{z}}^i d{\overline{z}}^j$$
where 
$$\alpha = \sum \alpha _{i{\overline j}}{d{z}}^i d{\overline{z}}^j$$
is the $(1,1)$ component of a form representing the class $\theta.$ This implies 
$${\operatorname Proj}[A^{(0)},A^{(0)}]_{1} = \omega (\theta _0, \theta _0).$$

\begin{remark} \label{rmk:canonical deformation} In \cite{NT1}, 4.8, we defined the canonical deformation of the trivial gerbe on a symplectic manifold. It is easy to see that the characteristic class $\theta$ of this deformation is equal to $\frac{1}{{\sqrt{-1}}\hbar}\omega.$ We see from Theorem \ref{thm:RW} that the first Rozansky-Witten class is an obstruction for the canonical stack deformation to be a sheaf of algebras.
\end{remark}
\section{Deformation complex of a stack as a DGLA}

In this section we will construct a DGLA whose Maurer-Cartan elements classify deformations of any stack (Theorem \ref{thm:deformations of stacks via dgla}). In order to that, we will start by noticing that a stack datum can be defined in terms of the simplicial nerve of a cover; if we replace the nerve by its first barycentric subdivision, we arrive at a notion of a descent datum for $\cL$ where $\cL$ is a simplicial sheaf of DGLAs (Definitions \ref{dfn:simplicial sheaf of dglas}, \ref{dfn:simplicial sheaf of dglas 1}). We reduce the problem to classifying such descent data in Proposition \ref{lemma:iso classes as L-stacks}. Then we replace our simplicial sheaf of DGLAs by a quasi-isomorphic acyclic simplicial sheaf of DGLAs. For the latter, classifying descent data is the same as classifying Maurer-Cartan elements of the DGLA of global sections, whence Theorem \ref{thm:deformations of stacks via dgla}. It states that deformations of a stack are classified by Maurer-Cartan elements of {\em De Rham-Sullivan forms with values in local Hochschild cochains of the twisted matrix algebra}.

\subsection{Twisted matrix algebras} \label{ss:Twisted matrix algebras}
 For any simplex $\sigma$ of the nerve of an open cover $M=\cup U_i$ corresponding to $U_{i_0}\cap \ldots \cap  U_{i_p}$, put $I_{\sigma}=\{i_o, \ldots, i_p\}$ and $U_{\sigma}=\cap_{i\in I}U_i.$ Define the algebra ${\operatorname {Matr}}^{\sigma}_{\operatorname{tw}}(\cA)$ whose elements are finite matrices
$$\sum_{i,j\in I_{\sigma}} a_{ij}E_{ij}$$ 
such that $a_{ij} \in \cA _i(U_{\sigma}). $ 
The product is defined by 
$$a_{ij}E_{ij}\cdot  a_{lk}E_{lk} = \delta_{jl} a_{ij}G_{ij}(a_{jk})c_{ijk}E_{ik}$$

We call a Hochschild $k$-cochain $D$ of ${\operatorname {Matr}}^{\sigma}_{\operatorname{tw}}(\cA)$ {\em local} if: 

a) For $k=0,$ $ D=\sum_{i\in I_{\sigma}} a_{i}E_{ii};$

b) for $k>0,$ $D(E_{i_1 j_1}, \ldots ,E_{i_k j_k}) = 0$
whenever $j_p \neq i_{p+1}$ for some $p$ between $1$ and $k-1;$

c) for $k>0,$ $D(E_{i_1 j_1}, \ldots ,E_{i_k j_k}) $ is a product of an element of $E_{i_1 j_k}$ and an element of $\cA$. 

Local cochains form a DGL subalgebra of all Hochschild cochains $C^{\bullet +1}({\operatorname {Matr}}^{\sigma}_{\operatorname{tw}}(\cA),\,{\operatorname {Matr}}^{\sigma}_{\operatorname{tw}}(\cA))$. Denote it by $\cL ^{H, \operatorname{local}}({\operatorname {Matr}}^{\sigma}_{\operatorname{tw}}(\cA)). $

\begin{remark} \label{rmk:local cochains and categories} 
It is easy to define a sheaf of categories on $U_{\sigma}$ whose complex of Hoschild cochains is exactly the complex of local Hochschild cochains above.
\end{remark}

\subsection{De Rham-Sullivan forms.} \label{De Rham-Sullivan forms}
 For any $p$-simplex $\sigma$ of the nerve of an open cover $M=\cup U_i$ corresponding to $U_{i_0}\cap \ldots \cap  U_{i_p}$, let 
$${\mathbb Q}[\Delta  _{\sigma}]={\mathbb Q}[t_{i_0}, \ldots, t_{i_p}]/(t_{i_0}+ \ldots + t_{i_p}-1)$$
and
$$\Omega ^{\bullet} [\Delta  _{\sigma}]={\mathbb Q}[t_{i_0}, \ldots, t_{i_p}]\{dt_{i_0}, \ldots, dt_{i_p}\}/(t_{i_0}+ \ldots + t_{i_p}-1,dt_{i_0}+ \ldots + dt_{i_p} )$$

As usual, given a sheaf $\cL$ on $M$, define {\em De Rham-Sullivan forms} with values in $\cL$ as collections $\omega _{\sigma}\in \Omega ^{\bullet} [\Delta  _{\sigma}]\otimes \cL(U_{\sigma})$ where $\sigma$ runs through all simplices, subject to $\omega _{\tau}|\Delta _{\sigma}=\omega _{\sigma}$ on $U_{\tau}$ whenever $\sigma \subset \tau$. De Rham-Sullivan forms form a complex with the differential $(\omega_{\sigma})\mapsto (d_{\operatorname{DR}}\omega_{\sigma}).$ We denote the space of all $k$-forms by $\ods ^k({\mathfrak U}, \cL)$, or simply by $\ods ^k({\mathfrak U})$ in the case when $\cL={\mathbb C}.$ The complex $(\ods ^{\bullet}({\mathfrak U}, \cL), d_{\operatorname{DR}})$ computes the \v{C}ech cohomology of $M$ with coefficients in $\cL$. Finally, put 
$$\ods ^k(M, \cL)=\limdir_{\mathfrak U}\ods ^k({\mathfrak U}, \cL)$$
where the limit is taken over the category of all open covers.

We need to say a few words about the functoriality of Hochschild cochains. Usually, given a morphism of algebras $A\to B$, there is no natural morphism between $C^{\bullet}(A,A)$ and $C^{\bullet}(B,B)$ (both map to $C^{\bullet}(A,B).$ Nevertheless, in our special case, there are maps ${\operatorname{Matr}}_{\operatorname{tw}}^{\sigma}\to {\operatorname{Matr}}_{\operatorname{tw}}^{\tau}$ on $U_{\tau}$ if $\sigma\subset\tau.$ These maps do induce morphisms of sheaves of {\em local} cochains on the open subset $U_{\tau}$ in the opposite direction; we call these morphisms {\em the restriction maps}. And, as before, we consider Hochschild cochain complexes already as sheaves of complexes. For example, in all the cases we are interested in, Hochschild cochains are given by multidifferential maps. 

\begin{definition}\label{dfn:DRS local}
Let $\odsb ({\mathfrak U}, \cL^{H,{\operatorname{local}}}(\mtwa))$ be the space of all collections 
$$D_{\sigma}\in \cL ^{H,{\operatorname{local}}}({\operatorname {Matr}}^{\sigma}_{\operatorname{tw}}(\cA))\otimes \Omega ^k(\Delta _{\sigma})$$
such that for $\sigma\subset\tau$ the restriction of the cochain $D_{\tau}|{\Delta _\sigma}$ to ${\operatorname {Matr}}^{\sigma}_{\operatorname{tw}}(\cA)$ is equal to $D_{\sigma}$ on $U_{\tau}.$ These spaces form a DGLA with the bracket $[(D_{\sigma}),\, (E_{\sigma})]=([D_{\sigma},\,E_{\sigma}])$ and the differential $(D_{\sigma})\mapsto ((d_{\operatorname{DR}}+\delta)D_{\sigma}).$ We put
$$\odsb (M, \cL^{H,{\operatorname{local}}}(\mtwa))=\limdir_{\mathfrak U}\odsb ({\mathfrak U}, \cL^{H,{\operatorname{local}}}(\mtwa))$$
\end{definition}

\begin{thm} \label {thm:deformations of stacks via dgla} Isomorphism classes of deformations of any stack $\cA$ are in one-to-one correspondence with isomorphism classes of Maurer-Cartan elements of the DGLA $\odsb (M, \cL^{H,{\operatorname{local}}}(\mtwa))$.
\end{thm}

The DGLAs above are examples of a structure that we call {\em a simplicial sheaf of DGLAs.}

\begin{definition}\label{dfn:simplicial sheaf of dglas} A simplicial sheaf $\cL$ is a collection of sheaves $\cL_{\sigma}$ on $U_{\sigma}$, together with morphisms of sheaves $r_{\sigma\tau}:\cL_{\tau}\to\cL_{\sigma}$ on $U_{\tau}$ for all $\sigma \subset  \tau$, such that  $r_{\sigma\tau}r_{\tau\theta}=r_{\sigma\theta}$ for any $\sigma \subset  \tau\subset  \theta.$ A simplicial sheaf of DGLAs $\cL$ is a simplicial sheaf such that all $\cL_{\sigma}$ are DGLAs and all $r_{\sigma\tau}$ are morphisms of DGLAs.
\end{definition}
\begin{definition}\label{dfn:simplicial sheaf of dglas 1}
For a simplicial sheaf of DGLAs $\cL$, {\em a descent datum} is a collection of Maurer-Cartan elements $\lambda _{\sigma} \in \hbar\cL ^1(U_{\sigma}[[\hbar]]),$ together with gauge transformations $G_{\sigma\tau}:r_{\sigma\tau}\lambda_{\tau}\to \lambda_{\sigma}$ on $U_{\tau}$ and two-morphisms $c_{\sigma\tau\theta}:G_{\sigma\tau}r_{\sigma\tau}(G_{\tau\theta})\to G_{\sigma\theta}$ on $U_{\theta}$ for any $\sigma \subset  \tau\subset  \theta,$ subject to
$$c_{\sigma\tau\omega}G_{\sigma\tau}(r_{\sigma\tau}(c_{\tau\theta\omega}))=c_{\sigma\theta\omega}c_{\sigma\tau\theta}$$
for any $\sigma \subset  \tau\subset  \theta\subset \omega.$
\end{definition}
We leave to the reader the definition of isomorphisms (and two-isomorphisms) of descent data. Given a simplicial sheaf $\cL$, and denoting the cover by ${\mathfrak U}$, one defines the cochain complex
$$C^{p}({\mathfrak U}, \cL)=\prod _{\sigma _{0}\subset \ldots \subset\sigma_{p}} \cL_{\sigma _0}(U_{\sigma_p})$$
Put
$$(d_0s)_{\sigma_{0}\ldots\sigma_{p+1}}=s _{\sigma_{1}\ldots\sigma_{p+1}};$$
$$(d_is)_{\sigma_{0}\ldots\sigma_{p+1}}=s _{\sigma_{0}\ldots{\widehat{\sigma _{i}}}\ldots\sigma_{p+1}},$$
$1\leq i \leq p;$
$$(d_{p+1}s)_{\sigma_{0}\ldots\sigma_{p+1}}=r_{\sigma_{p},\sigma_{p+1}} s_{\sigma_{0}\ldots\sigma_{p}}
$$
We leave to the reader the definition of the maps $s_i$. We see that $C^{\bullet}({\mathfrak U}, \cL)$ is a cosimplicial space. It is a cosimplicial DGLA if $\cL$ is a simplicial sheaf of DGLAs. 

Finally, note that, if a cover ${\mathfrak V}$ is a refinement of the cover ${\mathfrak U},$ then there is a morphism of cosimplicial spaces (DGLAs) 
$$C^{\bullet}({\mathfrak U}, \cL)\to C^{\bullet}({\mathfrak V}, \cL).$$
Let 
$$C^{\bullet}( \cL)=\limdir_{\mathfrak U}C^{\bullet}({\mathfrak U}, \cL).$$
We say that $\cL$ is {\em acyclic} if for every $q$ the cohomology of this complex is zero for $p>0$. 
\begin{definition}\label{dfn:Cech complex} The cochain complex $(C^{\bullet}({\mathfrak U}, \cL), \partial +d)$ where $\partial=\sum_{i=0}^n (-1)^id_i$ is called the \v{C}ech complex of $\cL$ with respect to the cover ${\mathfrak U}.$
\end{definition}
The collection of sheaves $\cL ^{H, \operatorname{local}}({\operatorname {Matr}}^{\sigma}_{\operatorname{tw}}(\cA)) $ forms a simplicial sheaf of DGLAs if one sets $r_{\sigma\tau}(\omega)$ to be the restriction of the $\omega $ to the algebra ${\operatorname {Matr}}^{\sigma}_{\operatorname{tw}}(\cA)$. We denote this simplicial sheaf of DGLAs by $\cL ^{H, \operatorname{local}}({\operatorname {Matr}}_{\operatorname{tw}}(\cA)) $.

\begin{proposition}\label{lemma:iso classes as L-stacks} Isomorphism classes of deformations over ${\mathfrak a}$ of any stack $\cA$ are in one-to-one correspondence with isomorphism classes of descent data of the Deligne two-groupoid of
$\cL ^{H, \operatorname{local}}({\operatorname {Matr}}_{\operatorname{tw}}(\cA))\otimes {\mathfrak a}.$
\end{proposition}

{\bf Proof}. Given a deformation, it defines a Maurer-Cartan element of $\cL ^{H, \operatorname{local}}({\operatorname {Matr}}_{\operatorname{tw}}^{\sigma}(\cA)) $ for every $\sigma$, namely the Hochschild cochain corresponding to the deformed product on ${\operatorname {Matr}}_{\operatorname{tw}}(\cA)$. It is immediate that this cochain is local. The restriction $r_{\sigma\tau}$ sends these cochains to each other, so a deformation of $\cA$ does define a descent datum for the Deligne two-groupoid of $\cL ^{H, \operatorname{local}}$. Conversely, to have such a descent datum is the same as to have a deformed stack datum ${\tilde \cA} _{\sigma}$ on every $U_{\sigma}$ (with respect to the cover by $U_i \cap U_{\sigma}= U_{\sigma},$ $i\in I_{\sigma}$), together with an isomorphism ${\tilde \cA} _{\tau}\to {\tilde \cA} _{\sigma}$ on $U_{\tau}$ for $\sigma \subset \tau$ and a two-isomorphism on $U_{\theta}$ for every $\sigma \subset \tau \subset \theta .$ But the cover consists of several copies of the same open set, which coincides with the entire space. All stack data with respect to such a cover are isomorphic to sheaves of rings; all stack isomorphisms are two-isomorphic to usual isomorphisms of sheaves. Trivializing the stacks ${\tilde \cA} _{\sigma}$ on $U_{\sigma}$ according to this, we see that isomorphism classes of such data are in one-to-one correspondence with isomorphism classes of the following:

1) a deformation ${\mathbb A}_{\sigma}$ of the sheaf of algebras $\cA _{i_0}$ on $U_{\sigma}$ where $I_\sigma = \{i_0, \ldots,i_p\}$;

2) an isomorphism of deformations $G_{\sigma \tau}:{\mathbb A}_{\tau}\to {\mathbb A}_{\sigma}|U_{\tau}$ for every $\sigma \subset \tau$;

3) an invertible element of $c_{\sigma \tau \rho} \in {\mathbb A}_{\sigma} (U_{\theta})$ for every $\sigma \subset \tau \subset \theta, $

satisfying the equations that we leave to the reader. Finally, 
one can establish a one-to-one correspondence between isomorphism classes of the above data and isomorphism classes of deformations of $\cA$. This is done using an explicit formula utilizing the fact that sequences $\sigma_0 \subset \ldots \subset \sigma_p$ are numbered by simplices of the barycentric subdivision of $\sigma _p$ (cf., for example, \cite{Seg}). More precisely, given a datum ${\mathbb A}_{\sigma}, G_{\sigma \tau}, c_{\sigma \tau \rho},$ we would like to construct a stack datum ${\mathbb A}_{i}, G_{ij}, c_{ijk}.$ We start by putting ${\mathbb A}_{i}={\mathbb A}_{(i)}$ and $G_{ij}=G_{(i),(ij)}G_{(j),(ij)}^{-1}.$ Now we want to guess a formula for $c_{ijk}.$ For that, observe that 
$$G_{(i),(ij)}={\operatorname{Ad}}(c_{(i), (ij), (ijk)})G_{(i), (ijk)}G_{(ij), (ijk)}^{-1}
$$
and 
$$G_{(j),(ij)}={\operatorname{Ad}}(c_{(j), (ij), (ijk)})G_{(j), (ijk)}G_{(ij), (ijk)}^{-1},
$$
therefore
$$G_{ij}={\operatorname{Ad}}(c_{(i), (ij), (ijk)})G_{(i), (ijk)}G_{(j), (ijk)}^{-1}{\operatorname{Ad}}(c_{(j), (ij), (ijk)}^{-1})
$$
We see that 
$$
G_{ij}G_{jk}={\operatorname{Ad}}(c_{ijk})G_{ik}$$
where
$$c_{ijk}=c_{(i), (ij), (ijk)}(G_{(i), (ijk)}G_{(j), (ijk)}^{-1})(c_{(j), (ij), (ijk)}^{-1}c_{(j), (jk), (ijk)})\times
$$
$$\times (G_{(i), (ijk)}G_{(k), (ijk)}^{-1})(c_{(k), (jk), (ijk)}^{-1}c_{(k), (ik), (ijk)})c_{(i), (ik), (ijk)}^{-1}
$$
(as one would expect, this is an alternated product of terms corresponding to the six faces of the first barycentric subdivision of the simplex $(ijk)$, in the natural order). One checks directly that the cocyclicity condition on the $c_{ijk}$'s holds. Furthermore, given an isomorphism $H_{\sigma}, b_{\sigma \tau}$ of the data ${\mathbb A}_{\sigma}, G_{\sigma \tau}, c_{\sigma \tau \rho}$ and ${\mathbb A}'_{\sigma}, G'_{\sigma \tau}, c'_{\sigma \tau \rho},$ one defines
$$H_i=H_{(i)}, \; b_{ij}=b_{(i),(ij)} b_{(j),(ij)}^{-1}$$
and checks that this is indeed an isomorphism of the corresponding data ${\mathbb A}_{i}, G_{ij}, c_{ijk}$ and ${\mathbb A}'_{i}, G'_{ij}, c'_{ijk}.$ This ends the proof of Lemma \ref{lemma:iso classes as L-stacks}.
\subsubsection{End of the proof of Theorem \ref{thm:deformations of stacks via dgla}} Define the simplicial sheaf of DGLAs as follows. Put 
$$\cL _{\sigma}= \cL ^{H,{\operatorname{local}}}({\operatorname {Matr}}^{\sigma}_{\operatorname{tw}}(\cA))\otimes \Omega ^k(\Delta _{\sigma}),$$
with the differential $d_{\operatorname{DR}}+\delta$ and transition homomorphisms 
$$r_{\sigma\tau}(D_{\tau})=D_{\tau}|\Delta _{\sigma}\; {\operatorname{restricted}}\;{\operatorname{to}} \;{\operatorname{Matr}}^{\sigma}_{\operatorname{tw}}(\cA).$$
We denote this simplicial sheaf of DGLAs by 
$${\underline {\Omega}} ^{\bullet}_{\operatorname{DRS}} (M, \cL^{H,{\operatorname{local}}}(\mtwa)).$$ 
It is acyclic as a simplicial sheaf. Therefore, by Proposition \ref{prop:quis of cosimplicial dglas}, isomorphism classes of descent data of its Deligne two-groupoid are in one-to-one correspondence with isomorphism classes of Maurer-Cartan elements of the DGLA $\odsb (M, \cL^{H,{\operatorname{local}}}(\mtwa))$, because the latter is its zero degree \v{C}ech cohomology. Now, the embedding
$$\cL^{H,{\operatorname{local}}}(\mtwa))\to {\underline {\Omega}} ^{\bullet}_{\operatorname{DRS}} (M, \cL^{H,{\operatorname{local}}}(\mtwa))$$
is a quasi-isomorphism of simplicial sheaves of DGLAs (the left hand side is the zero degree De Rham cohomology, and the higher De Rham cohomology vanishes locally). Again by Proposition \ref{prop:quis of cosimplicial dglas}, isomorphism classes of descent data are in one-to-one correspondence for the two simplicial sheaves of DGLAs above. 
\subsubsection{Another version of Theorem \ref{thm:deformations of stacks via dgla}} The language the previous subsection allows one to classify deformations of a given stack in terms of another DGLA which is a totalization of a cosimplicial DGLA. This is perhaps a little bit more consistent with the framework of \cite{H1}.
\begin{thm}\label{thm:deformations of stacks via dgla II} Isomorphism classes of deformations of a stack $\cA$ are in one-to-one correspondence with isomorphism classes of Maurer-Cartan elements of the DGLA $\Tot C^{\bullet}(\cL ^{H, \operatorname{local}}({\operatorname {Matr}}_{\operatorname{tw}}(\cA)))\otimes {\mathfrak a}.$
\end{thm}


\section{Deformations of a given gerbe} \label{s:Deformations of a given gerbe}
\subsection{} The aim of this section is to classify deformations of a given gerbe, trivial or not. As above, let $\cA$ be a gerbe on $M$; by $\cO _M$ we will denote the sheaf of smooth functions (in the $C^{\infty}$ case) or the holomorphic functions (in the complex analytic case).

The two-cocycle $c_{ijk}$ defining the gerbe belongs to the cohomology class in $H^2(M, {\cO}_M /2\pi i{\Bbb Z})$. Project this class onto $H^2(M, {\cO}_M / {\Bbb C})$. 
\begin{definition}\label{dfn:R}
We denote the above class in $H^2(M, {\cO}_M / {\Bbb C})$ by $R(\cA)$ or simply by $R$.
\end{definition}
The class $R$ can be represented by a two-form $R$ in $\ods^2(\cO_M/{\mathbb C})$, cf. \ref{De Rham-Sullivan forms}.
\begin{thm}\label{thm:classification of deformations of a gerbe}
Given a gerbe $\cA$ on a manifold $M$, the set of deformations over ${\mathfrak a}$ of $\cA$ up to isomorphism is in one-to-one correspondence with the set of equivalence classes of Maurer-Cartan elements of the DGLA $\ods^{\bullet}(M, C^{\bullet +1}(\cO_M, \cO_M))\otimes{\mathfrak m}$ with the differential $d_{\operatorname{DR}} + \delta + i_R.$
\end{thm}
Here $C^{\bullet +1}(\cO_M, \cO_M)$ is the sheaf of complexes of multi-differential Hochschild cochains of the jet algebra; $R \in \ods^2(M, \cO_M/{\mathbb C})$ is a form representing the class from Definition \ref{dfn:R}; $i_R$ is the Gerstenhaber bracket with the Hochschild zero-cochain $R$. Explicitly, if $R$ is an element of an algebra $A$, 
$$i_R D(a_1, \ldots, a_n)=\sum_{i=0}^n (-1)^i D(a_1, \ldots,a_i,\, R, \ldots, a_n).$$
In Theorem \ref{thm:classification of deformations of a gerbe} this operation is combined with the wedge multiplication on forms.

If the manifold $M$ is complex, we can formulate the theorem in terms of Dolbeault complexes, without resorting to De Rham-Sullivan forms.

\begin{thm}\label{thm:classification of deformations of a holomorphic gerbe}
Given a holomorphic gerbe $\cA$ on a complex manifold $M$, the set of deformations of $\cA$ over ${\mathfrak a}$ up to isomorphism is in one-to-one correspondence with the set of equivalence classes of Maurer-Cartan elements of the DGLA $A^{0, \bullet}(M, C^{\bullet +1}({\cO}_M, {\cO}_M))\otimes{\mathfrak m}$ with the differential ${\overline{\partial}} + \delta + i_R.$
\end{thm}
Here $R \in A^{0,2} (M, {\cO}_M/{\Bbb C})$ is a form representing the class from Definition \ref{dfn:R}; $i_R$ is the Gerstenhaber bracket with the Hochschild zero-cochain $R$.

 We start with a coordinate change that replaces twisted matrices by usual matrices, at a price of making the differential and the transition isomorphisms more complicated (Lemma \ref{lemma:first coordinate change}). The second coordinate change (\eqref{eq:second coordinate change} and up) allows to get rid of matrices altogether.

The rest of this section is devoted to the proof of the theorems above.

The plan of the proof is the following. Having reduced the problem of classifying deformations of a gerbe to the problem of classifying Maurer-Cartan elements of a DGLA (Theorem \ref{thm:deformations of stacks via dgla}), we will now simplify this DGLA.

\subsubsection{First coordinate change: untwisting the matrices} \label{ss:first coordinate change} Recall that we are working on a manifold M with an open cover $\{U_i\}_{i\in I}$ and a \v{C}ech two-cocycle $c_{ijk}$ with coefficients in $\cO^* _M$. 

In what follows, we will denote by $\Omega ^k(\Delta _{\sigma}, \cO (U_{\sigma}))$, etc. the space of forms on the simplex $\Delta _{\sigma}$ with values in $\cO (U_{\sigma})$, etc.

We start by observing that in the definition of De Rham-Sullivan forms one can replace {\em algebraic} $\cL$-valued forms $\Omega^{\bullet}(\Delta_{\sigma})\otimes {\cL}$ by {\em smooth} $\cL$-valued forms $\Omega^{\bullet}(\Delta_{\sigma},{\cL})$ where $\cL$ is the DGLA of local Hochschild cochains. Indeed, one DGLA embeds into the other quasi-isomorphically, and one can apply Proposition \ref{prop:quis of cosimplicial dglas}.

Locally, $c$ can be trivialized. Indeed, as in the proof of Proposition \ref{lemma:iso classes as L-stacks}, $c$ is a cocycle on $U_{\sigma}$ with respect to the cover of $U_{\sigma}$ by several copies of itself. We write
\begin{equation} \label{eq:first trivialization, 1}
c_{ijk} = h_{ij}(\sigma) h_{ik}(\sigma)^{-1} h_{jk}(\sigma)
\end{equation}
on $U_{\sigma}$ for a simplex $\sigma $, where  $h_{ij}$ are elements of $\Omega ^0(\Delta _{\sigma}, \cO (U_{\sigma}))$. As a consequence,
\begin{equation}\label{eq:first trivialization, 2}
d_{\operatorname{DR}}{\operatorname{log}} h_{ij}(\sigma)-d_{\operatorname{DR}} {\operatorname{log}} h_{ik}(\sigma)+d_{\operatorname{DR}} {\operatorname{log}} h_{jk}(\sigma)=0
\end{equation}

\begin{remark}\label{rmk:dependence} 
At this stage the cochains $h_{ij}(\sigma)$, $a_{i}(\sigma, \tau)$ can be chosen to be constant as functions on simplices. But later they will be required to satisfy Lemma \ref{lemma:choice of ai}, and for that they have to be dependent on the variables $t_i$.
\end{remark}

Note that two local trivializations of the two-cocycle $c$ differ by a one-cocycle which is itself locally trivial (by the same argument as the one before \eqref{eq:first trivialization, 1}). Therefore
\begin{equation} \label{eq:first trivialization, 3}
 h_{ij}(\sigma)  = a_{i}(\sigma, \tau) h_{ij}(\tau)a_{j}(\sigma, \tau)^{-1}
\end{equation}
on $ U_{\tau}$ where $a_i$ are some invertible elements of  $\Omega ^0(\Delta _{\sigma}, \cO (U_{\tau})).$ We have another local trivialization:
\begin{equation}\label{eq:first trivialization, 4}
d_{\operatorname{DR}} {\operatorname{log}} h_{ij}(\sigma) = \beta _i(\sigma)-\beta _j(\sigma)
\end{equation}
on $U_{\sigma}, $ where $\beta _i(\sigma)$ are elements of $\Omega ^1(\Delta _{\sigma}, \cO (U_{\sigma}))$. Now introduce the coordinate change
\begin{equation} \label{eq:first coordinate change}
 a_{ij}E_{ij} \mapsto a_{ij}h_{ij}(\sigma)E_{ij} 
\end{equation}
\begin{definition} \label{dfn:matrices}
By ${\operatorname {Matr}}_{\sigma}(\cA)$ we denote the sheaf on $U_{\sigma}$ whose elements are finite sums $\sum a_{ij}E_{ij}$ where $a_{ij}\in \cA_i.$ The multiplication is the usual matrix multiplication.
\end{definition}
One gets immediately 
\begin{lemma}\label{lemma:first coordinate change}
Put 
$$a(\sigma , \tau)={\operatorname {diag}}\,a_{i}(\sigma, \tau)$$
and
$$\beta(\sigma)={\operatorname {diag}}\,\beta_{i}(\sigma)$$

Consider the spaces of all collections 
$$D_{\sigma}\in \Omega ^k(\Delta _{\sigma}, \cL ^{H,{\operatorname{local}}}({\operatorname {Matr}}^{\sigma}(\cO)) )$$
such that for $\sigma\subset\tau$ the restriction of the cochain $D_{\tau}|{\sigma}$ to ${\operatorname {Matr}}^{\sigma}(\cA)$ is equal to ${\operatorname{Ad}}(a(\sigma,\tau))(D_{\sigma})$ on $U_{\tau}.$ These spaces form a DGLA with the bracket $[(D_{\sigma}),\, (E_{\sigma})]=([D_{\sigma},\,E_{\sigma}])$ and the differential $(D_{\sigma})\mapsto ((d_{\operatorname{DR}}+\delta+{\operatorname{ad}}(\beta(\sigma)))D_{\sigma}).$  The coordinate change \eqref{eq:first coordinate change} provides an isomorphism of this DGLA and the DGLA $\odsb (M, \cL(\mtwa))$ from Definition \ref{dfn:DRS local} (modified as in the beginning of \ref{ss:first coordinate change}).
\end{lemma}
\subsubsection{Second coordinate change} \label{ss:second coordinate change} We have succeeded in replacing the sheaf of DGLAs of Hochschild complexes of twisted matrices by the sheaf of DGLAs of Hochschild complexes of usual matrices, at a price of having more complicated differential and transition functions. Both involve conjugation (or commutator) with a diagonal matrix. Our next aim is to make these diagonal matrices have all the entries to be the same. This will allow us eventually to get rid of matrices altogether.

We already have one such diagonal matrix. Indeed, from \eqref{eq:first trivialization, 4} one concludes that 
\begin{equation}\label{eq:nabla beta is scalar}
d_{\operatorname{DR}} \beta _i(\sigma)=d_{\operatorname{DR}}  \beta _j(\sigma)
\end{equation}
and therefore
$$d_{\operatorname{DR}} \beta (\sigma) \in \Omega ^2(\Delta _{\sigma},\cO(U_{\sigma}) )$$
is well-defined. The other one is
\begin{equation} \label{eq:definition of gamma}
\gamma(\sigma, \tau)= d_{\operatorname{DR}} {\operatorname{log}}a_i(\sigma, \tau)-\beta _i (\sigma) + \beta _i(\tau)
\end{equation}
To see that this expression does not depend on $i$, apply $d_{\operatorname{DR}} {\operatorname{log}}$ to \eqref{eq:first trivialization, 3} and compare the result with \eqref{eq:first trivialization, 4}. Thus, we have a well-defined element
$$\gamma(\sigma, \tau) \in \Omega ^1(\Delta _{\sigma},\cO(U_{\tau}) ) .$$
Also, from \eqref{eq:first trivialization, 3} we observe that 
\begin{equation} \label{eq:definition of s}
s(\sigma, \tau, \theta)=a_i(\sigma, \tau)a_i(\sigma, \theta)^{-1}a_i(\tau, \theta)
\end{equation}
does not depend on $i$ and therefore defines an invertible element
$$s(\sigma , \tau , \theta) \in \Omega ^0(\Delta _{\sigma},\cO(U_{\theta}) ).$$
 The above cochains form a cocycle in the following sense:
\begin{equation}\label{eq:cech-nabla cocycle, 1}
d_{\operatorname{DR}} (d_{\operatorname{DR}} \beta)=0;
\end{equation}
\begin{equation}\label{eq:cech-nabla cocycle, 2}
d_{\operatorname{DR}} \beta(\sigma) - d_{\operatorname{DR}} \beta(\tau) = -d_{\operatorname{DR}} \gamma (\sigma, \tau);
\end{equation}
\begin{equation}\label{eq:cech-nabla cocycle, 3}
\gamma (\sigma , \tau) - \gamma (\sigma , \theta) + \gamma (\tau , \theta) = d_{\operatorname{DR}} {\operatorname {log}} s(\sigma , \tau , \theta);
\end{equation}
\begin{equation}\label{eq:cech-nabla cocycle, 4}
s(\sigma, \tau, \theta) s(\rho, \tau, \theta)^{-1}s(\rho, \sigma, \theta)s(\rho,\sigma, \tau )^{-1}=1
\end{equation}

\begin{lemma} \label{lemma:s,nablabeta,gamma}

The cohomology of the \v{C}ech bicomplex of the complex of simplicial sheaves  
$$\sigma\mapsto\Omega^0(\Delta_{\sigma}, \cO(U_{\sigma})) ^*\stackrel{d_{\operatorname{DR}}{\operatorname{log}}}{\longrightarrow} \Omega^1(\Delta_{\sigma}, \cO(U_{\sigma})) \stackrel{d_{\operatorname{DR}}}{\longrightarrow} \Omega^2(\Delta_{\sigma}, \cO(U_{\sigma})) \stackrel{d_{\operatorname{DR}}}{\longrightarrow}\ldots$$ 
is isomorphic to the \v{C}ech cohomology $H^{\bullet}(M, {\mathfrak U}; \cO^*_M)$ with respect to the cover ${\mathfrak U}.$ Under this isomorphism, the cohomology class of the cocycle $( d_{\operatorname{DR}} \beta, \gamma, s)$ of this complex becomes the cohomology class of the cocycle $c_{ijk}$.
\end{lemma}
The proof is straightforward, using the fact that sequences $\sigma_0 \subset \ldots \subset \sigma_p$ are numbered by simplices of the barycentric subdivision of $\sigma _p$ (cf. \cite{Seg}; compare with the proof of Proposition \ref{lemma:iso classes as L-stacks}) where a nonlinear version of the same argument is used).

From now on, we assume that the cover ${\mathfrak U}=\{U_i\}$ is good. We need another lemma to prooceed. 
\begin{lemma} \label{lemma:choice of ai}
The cochains $a_i (\sigma, \tau)$ can be chosen as follows:
$$a_i (\sigma, \tau)=a_0 (\sigma, \tau){\widetilde{a}}_i (\sigma, \tau)$$
where $a_0 (\sigma, \tau)$ does not depend on $i$ and ${\widetilde{a}}_i (\sigma, \tau)$ take values in the subgroup $\Omega ^0(\Delta _{\sigma}, {\C} \cdot 1)^*$.
\end{lemma}
{\bf Proof}. Choose local branches of the logarithm. We have from \eqref{eq:definition of s}
$$
{\operatorname {log}}a_i (\alpha, \sigma)-{\operatorname {log}}a_i (\alpha, \tau)+{\operatorname {log}}a_i ( \sigma , \tau)-{\operatorname {log}}s (\alpha, \sigma, \tau)=2\pi {\sqrt{-1}} N_i(\alpha, \sigma, \tau)
$$
where $N_i(\alpha, \sigma, \tau)$ are  constant integers. The \v{C}ech complex of the simplicial sheaf $\sigma \mapsto \Omega ^0(\Delta _{\sigma}, \cO _{U_{\sigma}}) $ is zero in positive degrees. Let $S$ be a contracting homotopy from this complex to its zero cohomology. Put
$$b_i (\sigma)=  {\operatorname{exp}}(S({\operatorname {log}}a_i (\alpha, \sigma)));$$
then
$$
b_i(\sigma)b_i(\tau)^{-1} = a_i(\sigma, \tau)^{-1} {\widetilde{a}}_i(\sigma, \tau) a(\sigma, \tau)
$$
where
$$
{\widetilde{a}}_i(\sigma, \tau)={\operatorname{exp}}(2\pi {\sqrt{-1}} S(N_i (\alpha, \sigma, \tau)))
$$
and
$$
a(\sigma, \tau)={\operatorname{exp}}(S(s(\alpha, \sigma, \tau)))
$$
Therefore we can, from the start, replace $h_{ij}(\sigma)$ by $b_i(\sigma)h_{ij}(\sigma)b_j(\sigma)^{-1}$ in \eqref{eq:first trivialization, 1}, and $a_i (\sigma, \tau)$ by $ {\widetilde{a}}_i(\sigma, \tau) a(\sigma, \tau)$ in \eqref{eq:first trivialization, 3}. This proves the lemma.

Now consider the operator 
$$i _{\beta(\sigma)}: {\Omega}^{\bullet} (\Delta _{\sigma}, C^{\bullet+1 }({\operatorname{Matr}}(\cO))) \to {\Omega}^{\bullet+1} (\Delta _{\sigma}, C^{\bullet}({\operatorname{Matr}}(\cO)))$$ 
This operator acts by the Gerstenhaber bracket (at the level of $C^{\bullet})$, combined with the wedge product at the level of $\Omega^{\bullet}$, with the cochain $\beta (\sigma) \in {\Omega}^{1} (\Delta _{\sigma}, C^0 ({\operatorname{Matr}}(\cO)))$. One has
$$[\delta, i _{\beta(\sigma)}] = {\operatorname{ad}} _{\beta(\sigma)}:{\Omega}^{\bullet} (\Delta _{\sigma}, C^{\bullet }({\operatorname{Matr}}(\cO))) \to {\Omega}^{\bullet+1} (\Delta _{\sigma}, C^{\bullet}({\operatorname{Matr}}(\cO)))$$
and 
$$[d_{\operatorname{DR}}, i _{\beta(\sigma)}] = i _{d_{\operatorname{DR}}\beta(\sigma)}$$ which is an operator 
$${\Omega}^{\bullet} (\Delta _{\sigma}, C^{\bullet+1 }({\operatorname{Matr}}(\cO))) \to {\Omega}^{\bullet+2} (\Delta _{\sigma}, C^{\bullet}({\operatorname{Matr}}(\cO)))$$
Now define the second coordinate change as
\begin{equation}\label{eq:second coordinate change}
{\operatorname{exp}}(i _{\beta(\sigma)})
\end{equation}
on ${\Omega}^{\bullet} (\Delta _{\sigma}, C^{\bullet}({\operatorname{Matr}}(\cO)))$. This coordinate change turns the DGLA from Lemma \ref{lemma:first coordinate change} into the following DGLA. Its elements are collections of elements
\begin{equation}\label{eq:new DGLA 1}
\omega _{\sigma}\in {\Omega}^{\bullet} (\Delta _{\sigma}, C^{\bullet}({\operatorname{Matr}}^{\sigma}(\cO ({U_{\sigma}}))))
\end{equation}
such that the restriction of $D_{\tau}|\Delta_{\sigma}$ to the subalgebra ${\operatorname{Matr}}^{\sigma}(\cO ({U_{\sigma}}))$ is equal to
\begin{equation}\label{eq:new DGLA 3}
 {\operatorname{exp}}(i _{\beta(\sigma)}-i_{\beta(\tau)}){\operatorname{Ad}}(a(\sigma, \tau))D_{\sigma};
\end{equation}
  the differential is
\begin{equation}\label{eq:new DGLA 2}
d_{\operatorname{DR}} + \delta +i _{d_{\operatorname{DR}} \beta(\sigma)}
\end{equation}
We can replace \eqref{eq:new DGLA 3} by
\begin{equation}\label{eq:new DGLA 4}
{\operatorname{exp}}(i _{\gamma(\sigma, \tau)}-i_{{d_{\operatorname{DR}}}{\operatorname{log}}{a_0}(\sigma, \tau)} - i_{{d_{\operatorname{DR}}}{\operatorname{log}}{\widetilde{a}}(\sigma, \tau)})){\operatorname{Ad}}(a_0(\sigma, \tau))D_{\sigma}
\end{equation}
where ${\widetilde{a}}(\sigma, \tau)={\operatorname{diag}}\,{\widetilde{a}}_i(\sigma, \tau)$ (cf. Lemma \ref{lemma:choice of ai}).
\subsection{Getting rid of matrices} Consider the morphism  
$$C^{\bullet}(\cO_{U_{\sigma}}) \to C^{\bullet}({\operatorname{Matr}}^{\sigma}(\cO_{U_{\sigma}})$$
defined as follows. Put ${\overline{\cO}}=\cO/{\Bbb C}$. Then for $D\in C^p(\cO,\cO),$ $D:{\overline{\cO}}^{\otimes p} \to {{\cO}},$ define
$${\widetilde{D}}(m_1a_1, \ldots, m_pa_p)= m_1\ldots m_p D(a_1, \ldots, a_p)$$
where $a_i \in \cO$ and $m_i \in M({\Bbb{C}}).$ The following is true:

a) the cochains ${\widetilde{D}}$ are invariant under isomorphisms ${\operatorname{Ad}}(m)$ for $m\in GL({\Bbb C})$; 

b) the cochains ${\widetilde{D}}$ become zero after substituting and argument from $M({\Bbb C})$.

It is well known that the map $D\mapsto {\widetilde{D}}$ is a quasi-isomorphism with respect to the Hochschild differential $\delta$. Therefore this map establishes a quasi-isomorphism of the DGLA from \eqref{eq:new DGLA 1}, \eqref{eq:new DGLA 2}, \eqref{eq:new DGLA 3},  \eqref{eq:new DGLA 4} with the following DGLA: its elements are collections $D_{\sigma}\in \Omega ^{\bullet} (\Delta _{\sigma}, C^{\bullet +1}(\cO ({U_{\sigma}})))$ such that 
\begin{equation}\label{eq:neznayu 2}
D_{\tau}|{\Delta _{\sigma}}={\operatorname{exp}}(i_{\gamma (\sigma, \tau)} - i_{d{\operatorname{log}}\,{a_0}(\sigma, \tau)})D_{\sigma}
\end{equation}
on $U_{\tau}$, 
with the differential
\begin{equation} \label{eq:neznayu 1}
d_{\operatorname{DR}} + \delta + i_{ d_{\operatorname{DR}}\beta(\sigma)}.
\end{equation}

Now consider any cocycle $r(\sigma)\in {\Omega }^2 (U_{\sigma}, \cO/{\Bbb C})$, $t(\sigma, \tau)\in \Omega ^1 ( U_{\tau}, \cO/{\Bbb C});$
$$r(\sigma)-r(\tau) + t(\sigma, \tau) = 0;$$
$$t(\sigma, \tau)-t(\sigma, \theta)+t(\tau, \theta)=0$$
Such a cocycle defines a of DGLA of collections $D_{\sigma}$ as above, where \eqref{eq:neznayu 2} gets replaced by 
\begin{equation}\label{eq:neznayu 21}
D_{\tau}|{\Delta _{\sigma}}={\operatorname{exp}}(i_{t(\sigma, \tau)})D_{\sigma}
\end{equation}
and the differential is $d_{\operatorname{DR}} + \delta + i_{r(\sigma)}$  If two cocycles differ by the differential of $u(\sigma) \in \Omega^1(\Delta ^{\sigma}, \cO(U_{\sigma})/{\Bbb  C})$, then operators ${\operatorname {exp}}(i_{u(\sigma)})$ define an isomorphism of DGLAs. Finally, put $r(\sigma)= \beta(\sigma)$ and $t(\sigma, \tau)=\gamma (\sigma, \tau)-d{\operatorname{log}}{\,a_0}(\sigma, \tau)$. This is a cocycle of \v{C}$^{\bullet} (M, {\cA}_M(\cO/{\Bbb C}))$. It lies in the cohomology class of the cocycle $({\operatorname{log}}\,s, \; \gamma, \;d_{\operatorname{DR}} \beta)$ from Lemma \ref{lemma:s,nablabeta,gamma}. Now replace this cocycle by a cohomologous cocycle which has $t=0$.

This proves that isomorphism classes of deformations of a gerbe $\cA$ are in one-to-one correspondence with isomorphism classes of Maurer-Cartan elements of the DGLA of collections of cochains 
$$D_{\sigma}\in \Omega^{\bullet}(\Delta _{\sigma}, C^{\bullet +1}(\cO_{U_{\sigma}}, \cO_{U_{\sigma}}))$$
such that $D_{\sigma}|U_{\tau}=D_{\tau};$ the differential is $d_{\operatorname{DR}}+\delta+i_R$ where $R\in \Omega ^2 _{\operatorname{DRS}}(M,\cO/{\mathbb C})$ represents the class $R$ as defined in the beginning of this section. To pass to the DGLA of Dolbeault forms (Theorem \ref{thm:classification of deformations of a holomorphic gerbe}), we apply Proposition \ref{prop:quis of cosimplicial dglas}.
\subsubsection{The jet formulation} Theorem \ref{thm:classification of deformations of a gerbe} also admits a formulation in the language of jets. As above, let $J_M$ be the bundle of algebras whose fiber at a point is the algebra of jets of $C^{\infty}$, resp. holomorphic, functions on $M$ at this point; this bundle has the canonical flat connection $\nc$. Horizontal sections of $J_M$ correspond to smooth, resp. holomorphic, functions. 

The two-cocycle $c_{ijk}$ defining the gerbe belongs to the cohomology class in $H^2(M, {\cO}_M /2\pi i{\Bbb Z})$. Project this class onto $H^2(M, {\cO}_M / {\Bbb C})$ and denote the result by $R$ (as in Definition \ref{dfn:R}).
The class $R$ can be represented by a two-form $R$ in $A^2 (M, J_M/{\Bbb C})$.
\begin{thm}\label{thm:classification of deformations of a gerbe, jet version}
Given a gerbe $\cA$ on a manifold $M$, the set of deformations of $\cA$ over ${\mathfrak a}$ up to isomorphism is in one-to-one correspondence with the set of equivalence classes of Maurer-Cartan elements of the DGLA $A^{\bullet}(M, C^{\bullet +1}(J_M, J_M))\otimes{\mathfrak m}$ with the differential $\nc + \delta + i_R.$
\end{thm}
Here $C^{\bullet +1}(J_M, J_M)$ is the complex of vector bundles of Hochschild cochains of the jet algebra; $R \in A^2 (M, J_M/{\Bbb C})$ is a form representing the class from Definition \ref{dfn:R}; $i_R$ is the Gerstenhaber bracket with the Hochschild zero-cochain $R$. The proof follows from a simple application of Proposition \ref{prop:quis of cosimplicial dglas}.
\section{Deformations of gerbes on symplectic manifolds}\label{s:Deformations of gerbes on symplectic manifolds}
\subsection{} For a gerbe on $M$ defined by a cocycle $c$, we denote by $c$ the class of this cocycle in $H^2(M, {\cO} _M /2\pi i{\Bbb Z})$ and by $\partial c$ its boundary in $H^3(M, 2\pi i{\Bbb Z})$.
\begin{thm}\label{thm:symplectic classification} Let $\cA$ be a gerbe on a symplectic manifold $(M, \omega)$. The set of isomorphism classes of deformations of $\cA$ compatible to $\omega$:

a) is empty if the image of the class $\partial c$ under the map $H^3(M, 2\pi i{\Bbb Z})\to H^3(M, {\Bbb C})$ is non-zero;

b) is in one-to-one correspondence with the space $Def(M,\omega)$ (Theorem \ref{thm:classification of deformations of the trivial gerbe, symplectic case})
 if the image of the class $\partial c$ under the map $H^3(M, 2\pi i{\Bbb Z})\to H^3(M, {\Bbb C})$ is zero.
\end{thm}
Let $R$ be the projection of $c$ to $H^2(M, {\cO} _M /{\Bbb C})$, as in Definition \ref{dfn:R}.
\begin{thm}\label{thm:symplectic holomorphic classification} Let $\cA$ be a gerbe on a complex symplectic manifold $(M, \omega)$. The set of isomorphism classes of deformations of $\cA$ compatible to $\omega$ is:

a) is empty if $R\neq 0$;

b) is in one-to-one correspondence with the space $Def(M,\omega)$ 
 if $R=0$.
\end{thm}
{\bf Proof.} The arguments from the proof of Theorem \ref{thm:classification of deformations of the trivial gerbe, symplectic case} show that deformations of a gerbe are classified exactly as in \eqref{eq:Fedosov 1}-\eqref{eq:Fedosov 4}, with one exception: equation \eqref{eq:Fedosov 2} should be replaced by the requirement that the class of $c$ modulo $A^2 (M, \C + \hbar{\operatorname{gr}} J)[[\hbar]]$ should coincide with $R$
where $R$ is a form defined before Theorem \ref{thm:classification of deformations of a gerbe}. Therefore, if $R=0,$ the classification goes unchanged; if $R\neq 0$ in $H^2(M, {\cO} _M /{\Bbb C})$, then 
\begin{equation}\label{eq:flatness mod h with R}
\no A^{(0)}+{\frac{1}{2}}[A^{(-1)},A^{(-1)}]_{2}=R 
\end{equation}
shows that no connection $\nabla$ exists.

\end{document}